%%%%%%%%%%%%%%%%%%%%%%%%%%%%%%%%%%%%%%%%%%%%%%%%%%%%%%
% The moduli space of $(1,11)$-polarized abelian surfaces is 
%             unirational
%
% by Mark Gross & Sorin Popescu
% 02/01/98
% ( Plain TeX with diagrams.tex)
%%%%%%%%%%%%%%%%%%%%%%%%%%%%%%%%%%%%%%%%%%%%%%%%%%%%%%
\font\hd=cmbx10 scaled\magstep1
\magnification=\magstep 1
\input diagrams.tex

% defs2.tex
%================================================
% F o n t s
%================================================
% define \Cal as control sequence with an argument for \fam2 in math mode

%------------------------------------------------
% For the following definitions of font families and
% corresponding control sequences see:
% AMSFonts 2.0 User's Guide page 12.
% There must no other families be defined before these,
% if you want to refer to the \fam numbers
% of the families starting from number 8 for \msamfam.
%------------------------------------------------
% define \Msa as control sequence with an argument for msam in math mode:
\font\tenmsam=msam10
\font\sevenmsam=msam7
\font\fivemsam=msam5
\newfam\msamfam         % family "8 (if no other \newfam before)
\textfont\msamfam\tenmsam
\scriptfont\msamfam\sevenmsam
\scriptscriptfont\msamfam\fivemsam

%
% define \MSAM as control sequence to switch to the msam font.

%------------------------------------------------
% define \Bbb#1 as control sequence with an argument for msbm in math mode:
\font\tenmsbm=msbm10
\font\sevenmsbm=msbm7
\font\fivemsbm=msbm5
\newfam\msbmfam         % family "9 (if no other \newfam before)
\textfont\msbmfam\tenmsbm
\scriptfont\msbmfam\sevenmsbm
\scriptscriptfont\msbmfam\fivemsbm

%
% define \MSBM as control sequence to switch to the msbm font.

%------------------------------------------------
% define \frak#1 as control sequence with an argument for eufm in math mode:
\font\teneufm=eufm10
\font\seveneufm=eufm7
\font\fiveeufm=eufm5
\newfam\eufmfam        % family "A (if no other \newfam before)
\textfont\eufmfam\teneufm
\scriptfont\eufmfam\seveneufm
\scriptscriptfont\eufmfam\fiveeufm

%
% define \Fraktur as control sequence to switch to the eufm font.

%------------------------------------------------
%
%===========================================================
%     D e f i n i n g    N e w   S y m b o l s
%===========================================================
% Definition of new symbols from the AMS-Fonts
% by the AMS-identification number (in hex):
%
% The first digit is the AMS family number:
%  1 for the MSAM family and 2 for the MSBM family.
%
% The second digit is the class number (for
% example 2 for a binary operator, or 0 for
% an ordinary character). Same as in Plain TEX, but in Plain TEX
% this is the first digit, and the second is the family number.
%
% The last 2 digits give the position in the font table.
%
% Therefore the macro must replace the first digit by the corresponding
% Plain TEX font family number and then interchange the first two digits
% in order to switch to the Plain Tex identification (ID-) number. 
% (See Schwarz, Einf?hrung in TEX page 256 \mathchar (or The TEX Book
% p. 154)  and page 264 \newfam.
%
%================================================
%
% The following macro "\newsymbol " works for MSAM, MSBM und EUFM 
% automatically using the AMS-Font ID-number 1 for MSAM, 2 for MSBS,
% 3 for EUFM similar to the AMS-Macro \newsymbol
%-------------------------------------------------------------
\newcount\amsfamcount % define new registers to
\newcount\classcount   % store the digits for conversion 
\newcount\positioncount
\newcount\codecount
\newcount\n             % store the #2 parameter
%-------------------------
\def\newsymbol#1#2#3#4#5{               % #1 is the name of the new
%                                       % macro to be defined
\n="#2                                  % #2 is the AMS font family number
\ifnum\n=1 \amsfamcount=\msamfam\else   % translate it to the correct
\ifnum\n=2 \amsfamcount=\msbmfam\else   % PlainTex font family number
\ifnum\n=3 \amsfamcount=\eufmfam
\fi\fi\fi
\multiply\amsfamcount by "100           % make it the third hex digit 
\classcount="#3                 % #3 is the class number (same as in Plain)
\multiply\classcount by "1000           % make it the fourth hex digit
\positioncount="#4#5            % #4#5 is the hex position in the fonttable
\codecount=\classcount                  % add the three numbers to obtain 
\advance\codecount by \amsfamcount      % the Plain TEX identification number
\advance\codecount by \positioncount
\mathchardef#1=\codecount}              % define #1 to become a character
%                                       % with this identification number
%========================================================================
% The following macro " \newmathsymbol " takes as
%    #1  the name of the new control sequence,
%    #2  the Plain TEX number of the font family, a decimal number 0-{15} or 
%        in the invariant form \msamfam, \msbmfam ,\eufmfam, ... which is 
%        independent of the sequence in which those families are loaded
%    #3  the Plain TEX class number as hex digit 0-7
%  #4#5  the position number of the character in the font table as a
%        two digit hex number 00-7F.
%------------------------------------------------------------------------
\newcount\famcnt % define new registers to
\newcount\classcnt   % store the digits for conversion %
\newcount\positioncnt
\newcount\codecnt
%-------------------------
\def\newmathsymbol#1#2#3#4#5{          % #1 is the name of the new
%                               % macro to be defined
\famcnt=#2                      % #2 is the font family number in decimal
\multiply\famcnt by "100        % make it the third hex digit 
\classcnt="#3                   % #3 is the class number (same as in Plain)
\multiply\classcnt by "1000     % make it the fourth hex digit
\positioncnt="#4#5              % #4#5 is the hex position in the fonttable
\codecnt=\classcnt              % add the three numbers to obtain 
\advance\codecnt by \famcnt     % the Plain TEX identification number
\advance\codecnt by \positioncnt
\mathchardef#1=\codecnt}        % define #1 to become a character
%                               % with this identification number
%=======================================================================

%rational short map
\def\rto{\raise.5ex\hbox{$\scriptscriptstyle ---\!\!\!>$}}

\def\Pfour{{\bf P^4}}
\def\Pthree{{\bf P^3}}

\def\Pone{{\bf P^1}}
\def\P{{\bf P}}

\def\A{{\cal A}}

\def\O{{\cal O}}
\def\L{{\cal L}}

\def\H{{\cal H}}
\def\G{{\cal G}}

\def\I{{\cal I}}

\def\E{{\cal E}}
\def\K{{\cal K}}

\def\HHH{{\bf H}}

\def\EY{\E(1)|_Y}

\def\PEY{\P(\EY)}

\def\OP{\O_{\PEY}(1)}
\def\OP1{\O_{\Pone}}

\def\im{\mathop{\rm im}}

\def\char{\mathop{\rm char}}
\def\Gr{\mathop{\rm Gr}\nolimits}

\def\mod{\mathop{\rm mod}}

\def\hom{\mathop{\rm Hom}\nolimits}

\def\rank{\mathop{\rm rank}\nolimits}

\def\boldz{{\bf Z}}

\def\dual#1{{#1}^{\scriptscriptstyle \vee}}
\newsymbol\SEMI226F
\def\rtimes{\mathop{\SEMI}}

%\def\normalbaselines{\baselineskip20pt
%\lineskip3pt \lineskiplimit3pt }

% defs3.tex
%================================================
%   M a t h. O p e r a t o r s
%================================================
%

%

%

%

%

%
\def\ker{\mathop{\rm ker}}
\def\im{\mathop{\rm im}}
\def\dim{\mathop{\rm dim}}
%

%
%==========================================================
%        O t h e r  C o n t r o l   S e q u e n c e s
%==========================================================
% qed in textmode
\newsymbol\QED1003
\def\qed{\hbox{$\QED$}}
%----------------------------------------------------------
% semi-direct product
\newsymbol\SEMI226F
\def\rtimes{\mathop{\SEMI}}
%------------------------------------------------------------
% define font "line10" as "lineten" and generate a control
% sequence "arrow" which prints the wanted character from
%line10 by its octal number:
\font\lineten=line10
\def\arrow#1{\lineten\char'#1}
%------------------------------------------------------------
% bold >= from AMS-TeX

%------------------------------------------------------------
% rightarrow with symbols above: #1 symbols on arrow:

%------------------------------------------------
% leftarrow with symbols above:

%----------------------------------------------------------------
% downarrows with symbols :

%-------------------------------------------------------------------
% upwards arrows with symbols :

%--------------------------------------------------------------
% hooked rightarrow with symbols :

%--------------------------------------------------------------------
% long rightarrow with symbols above:
% #1: symbols on arrow,  #2: length of the extension in cm.
\def\harrowext#1{\hskip-3.5pt\raise2.3pt\hbox
to#1{\hrulefill}\hskip-3.5pt}
%

%---------------------------------------------------------------
%long left with symbols above:

%----------------------------------------------------------------
% nwarrow with symbols above:

%----------------------------------------------------------------
% nwarrow with symbols above:

%----------------------------------------------------------------
% long north west arrow
%
\def\longnwarrow{\hskip-15pt\lower6pt\hbox{\vbox{\hbox{$\nwarrow$}
\vskip-1.5pt\hbox{\hskip7.95pt\arrow100 }}}}
%----------------------------------------------------------------
%rational short map
\def\rto{\raise.5ex\hbox{$\scriptscriptstyle ---\!\!\!>$}}
%-----------------------------------------------------------------
% New line in hmode:

%
% newline more general

%
% uparrow two characters long
\def\twoup{\vbox{\hbox{\hskip-2pt$\uparrow$}
\vskip-3pt\hbox{\hskip-1pt$|$}}}
%
% uparrow three characters long

%
% shortmid from AMS-TeX
\newsymbol\shortmid2370
%==========================================================
%

%============================================================
%       S h e a v e s
%============================================================
% ideal sheaf:

%
% twisted structural sheaf:

%
%  structural sheaf:

%
% twisted structural sheaf (simplified) :

%
% sheaf of differentials (simplified) :

%

%%%%%%%%%%%%%%%%%%%%

\def\QQ{{\bf Q}}

\def\vt{\vartheta}
\def\dim{{\rm dim}\,}

\def\Pf{{\rm Pf}\,}

% Macros from defs3.tex 
\def \PSL {\mathop{\rm PSL}\nolimits}
\def \SL {\mathop{\rm SL}\nolimits}
\def \GL {\mathop{\rm GL}\nolimits}
\def\PGL {\mathop{\rm PGL}\nolimits}
\def\vt{\vartheta}
\def\dim{\mathop{\rm dim}\nolimits}

\def\Pf{{\rm Pf}\,}
\def \Sec {\mathop{\rm Sec}\nolimits}

\font\abst=cmr9

\font\smallrm=cmcsc8

\headline={\ifodd\pageno \ifnum\pageno>1 \smallrm \hfil 
The moduli space of $(1,11)$-polarized abelian surfaces is 
unirational               % running head for right-hand page is title in caps
\hfil\folio \else\hfill\fi \else \smallrm \folio \hfill
M.~Gross and S.~Popescu   % running head for left-hand page is authors in caps
\hfill\fi} \footline={\hss}   % footline is blank

% begin.tex

%=========================================================================%
% load begin.tex only once, but keep count to match \bye commands
%=========================================================================%

\ifx\begin\undefined\else\global\advance\srcdepth by
1\expandafter \fi

\def\begin{}
\newcount\srcdepth
\srcdepth=1

\outer\def\bye{\global\advance\srcdepth by -1
  \ifnum\srcdepth=0
    \def\endcmd{\vfill\eject\nopagenumbers\par\vfill\supereject\end}
  \else\def\endcmd{}\fi
  \endcmd
}

%=========================================================================%
% initialize TeX
%=========================================================================%

\baselineskip=15pt
%\hsize = 6truein
%\hoffset = 0.5truein
%\vsize = 8.5truein
\voffset = 0.2truein
\emergencystretch = 0.05\hsize
\parskip=3pt plus1pt minus.5pt
\overfullrule=0pt

\newif\ifblackboardbold

% comment out the following line if AMS msbm fonts aren't available
\blackboardboldtrue

%=========================================================================%
% select fonts
%=========================================================================%

\font\sectionfont=cmbx10 scaled\magstephalf

% Establish AMS blackboard bold fonts without using amssym.def, amssym.tex

\newfam\bboldfam
\ifblackboardbold
\font\tenbbold=msbm10
\font\sevenbbold=msbm7
\font\fivebbold=msbm5
\textfont\bboldfam=\tenbbold
\scriptfont\bboldfam=\sevenbbold
\scriptscriptfont\bboldfam=\fivebbold
\def\bbold{\fam\bboldfam\tenbbold}
\else
\def\bbold{\bf}
\fi

%=========================================================================%
% font size-changing command ("A Beginner's Book of TeX" p35, p275)
%=========================================================================%

\font\Arm=cmr8
\font\Ai=cmmi8
\font\Asy=cmsy8
\font\Abf=cmbx8
\font\Brm=cmr6
\font\Bi=cmmi6
\font\Bsy=cmsy6
\font\Bbf=cmbx6
\font\Crm=cmr5
\font\Ci=cmmi5
\font\Csy=cmsy5
\font\Cbf=cmbx5

\ifblackboardbold
\font\Abbold=msbm10 at 8pt
\font\Bbbold=msbm7 at 6pt
\font\Cbbold=msbm5
\fi

\def\smallmath{%
\textfont0=\Arm \scriptfont0=\Brm \scriptscriptfont0=\Crm
\textfont1=\Ai \scriptfont1=\Bi \scriptscriptfont1=\Ci
\textfont2=\Asy \scriptfont2=\Bsy \scriptscriptfont2=\Csy
\textfont\bffam=\Abf \scriptfont\bffam=\Bbf \scriptscriptfont\bffam=\Cbf
\def\rm{\fam0\Arm}\def\mit{\fam1}\def\oldstyle{\fam1\Ai}%
\def\bf{\fam\bffam\Abf}%
\ifblackboardbold
\textfont\bboldfam=\Abbold
\scriptfont\bboldfam=\Bbbold
\scriptscriptfont\bboldfam=\Cbbold
\def\bbold{\fam\bboldfam\Abbold}%
\fi
}

%=========================================================================%
% single-pass symbolic theorem labeling
%=========================================================================%

% Because this is a single-pass mechanism with no .aux file, forward
% references need to be declared in advance:

%   \forward{thm:main}{Theorem}{1.1}

% This is also the mechanism for "timely" declaration of labels, which
% will usually be buried within the corresponding theorem macros.
% A warning is issued if a label redeclaration is inconsistent, allowing
% forward references to be manually fixed.

%   \ref{thm:main} produces "Theorem~1.1"
%   \refs{thm:main} produces "Theorems~1.1"
%   \refn{thm:main} produces "1.1"

% Some TeX adapted from "The Advanced TeXbook" by David Salomon, chapter 9.

% Implementers: The code for \forward is subtle. Its second argument must
% be provided literally, e.g. "Theorem" rather that "\capitalize{theorem}".
% Its third argument must either be literal or a macro that expands
% directly to a literal, e.g. "\edef\numtoks{\number\proccount}".
% This use of \edef cannot be replaced by \def, which defers expansion.
% Failure to follow these rules will cause spurious warnings that forward
% references are inconsistent, when they are in fact consistent after
% expansion. Note the "Towers of Palo Alto" recreational math problem
% involving the iterated use of \expandafter to expand the first argument
% to \forwardsub before calling it.

\newlinechar=`@
\def\forwardmsg#1#2#3{\immediate\write16{@*!*!*!* forward reference should
be: @\noexpand\forward{#1}{#2}{#3}@}}
\def\nodefmsg#1{\immediate\write16{@*!*!*!* #1 is an undefined reference@}}

\def\forwardsub#1#2{\def\newref{{#2}{#1}}}

\def\forward#1#2#3{%
\expandafter\expandafter\expandafter\forwardsub\expandafter{#3}{#2}
\expandafter\ifx\csname#1\endcsname\relax\else%
\expandafter\ifx\csname#1\endcsname\newref\else%
\forwardmsg{#1}{#2}{#3}\fi\fi%
\expandafter\let\csname#1\endcsname\newref}

\def\firstarg#1{\expandafter\argone #1}\def\argone#1#2{#1}
\def\secondarg#1{\expandafter\argtwo #1}\def\argtwo#1#2{#2}

\def\ref#1{\expandafter\ifx\csname#1\endcsname\relax
  {\nodefmsg{#1}\bf`#1'}\else
  \expandafter\firstarg\csname#1\endcsname
  ~\expandafter\secondarg\csname#1\endcsname\fi}

\def\refs#1{\expandafter\ifx\csname#1\endcsname\relax
  {\nodefmsg{#1}\bf`#1'}\else
  \expandafter\firstarg\csname #1\endcsname
  s~\expandafter\secondarg\csname#1\endcsname\fi}

\def\refn#1{\expandafter\ifx\csname#1\endcsname\relax
  {\nodefmsg{#1}\bf`#1'}\else
  \expandafter\secondarg\csname #1\endcsname\fi}

%=========================================================================%
% widow control
%=========================================================================%

% usage:
% \widow{.2} % start new page if <.2 page left

\def\widow#1{\vskip 0pt plus#1\vsize\goodbreak\vskip 0pt plus-#1\vsize}

%=========================================================================%
% sections and theorems
%=========================================================================%

% use \showlabels or \showlabelsabove to display section and theorem labels

\def\marginlabel#1{}

\def\showlabelsabove{
\font\labelfont=cmss10 at 6pt
\def\marginlabel##1{\rlap{\smash{\raise 10pt\hbox{\labelfont##1}}}}
}

\newcount\seccount
\newcount\proccount
\seccount=0
\proccount=0

\def\stdskip{\vskip 9pt plus3pt minus 3pt}
\def\stdbreak{\par\removelastskip\penalty-100\stdskip}

\def\proof{\stdbreak\noindent{\it Proof}}

\def\qed{\vrule height 1.2ex width .9ex depth .1ex}

\def\Box{
  \ifmmode\eqno\qed
  \else\ifvmode\removelastskip\line{\hfil\qed}
  \else\unskip\quad\hskip-\hsize
    \hbox{}\hskip\hsize minus 1em\qed\par
  \fi\stdbreak\fi}

\def\references{
  \removelastskip
  \widow{.05}
  \vskip 24pt plus 6pt minus 6 pt
  \leftline{\sectionfont References}
  \nobreak\stdskip\noindent}

\def\ifempty#1#2\endB{\ifx#1\endA}
\def\makeref#1#2#3{\ifempty#1\endA\endB\else\forward{#1}{#2}{#3}\fi}

\outer\def\section#1 #2\par{
  \removelastskip
  \global\advance\seccount by 1
  \global\proccount=0\relax
                \edef\numtoks{\number\seccount}
  \makeref{#1}{Section}{\numtoks}
  \widow{.05}
  \vskip 24pt plus 6pt minus 6 pt
  \message{#2}
  \leftline{\marginlabel{#1}\sectionfont\numtoks\quad #2}
  \nobreak\stdskip}

\def\proclamation#1#2{
  \outer\expandafter\def\csname#1\endcsname##1 ##2\par{
  \stdbreak
  \advance\proccount by 1
  \edef\numtoks{\number\seccount.\number\proccount}
  \makeref{##1}{#2}{\numtoks}
  \noindent{\marginlabel{##1}\bf #2 \numtoks\enspace}
  {\sl##2\par}
  \stdbreak}}

\def\othernumbered#1#2{
  \outer\expandafter\def\csname#1\endcsname##1{
  \stdbreak
  \advance\proccount by 1
  \edef\numtoks{\number\seccount.\number\proccount}
  \makeref{##1}{#2}{\numtoks}
  \noindent{\marginlabel{##1}\bf #2 \numtoks\enspace}}}

\proclamation{definition}{Definition}
\proclamation{lemma}{Lemma}
\proclamation{proposition}{Proposition}
\proclamation{theorem}{Theorem}
\proclamation{corollary}{Corollary}
\proclamation{conjecture}{Conjecture}

\othernumbered{example}{Example}
\othernumbered{remark}{Remark}
\othernumbered{construction}{Construction}
\othernumbered{claim}{Claim}
\othernumbered{question}{Question}

%=========================================================================%
% enable postscript illustrations using epsf.tex
%=========================================================================%

% Usage:
% \draw{70}{fig}{} % draw fig.eps at 70% scale
% \draw{999}{fig}{} % draw fig.eps scaled to width of page

% Optional third argument can be multiple calls to \figtext; see below.
% More generally, the third argument is read in vertical mode, with the
% reference point at the lower left corner of the eps picture, whose
% dimensions are contained in the dimen registers \drawx and \drawy.
% This enables using TeX to generate the text that goes with the picture.
% To request that the picture be widened to respect the added text, 
% examine and modify the dimen registers \ngap, \egap, \sgap, \wgap.
% This is done automatically by the \figtext macro.

% These macros rely on "epsf.tex" which is the lowest level interface
% available for including encapsulated Postscript files in TeX documents.
% Rather that manually reading the .eps file to compute the nominal size,
% the \epsfbox macro is called twice, and two of its internal registers
% are examined after the first call. A major change to epsf.tex (unlikely)
% will require changes here. 

%\input epsf

\newcount\figcount
\figcount=0
\newbox\drawing
\newcount\drawbp
\newdimen\drawx
\newdimen\drawy
\newdimen\ngap
\newdimen\sgap
\newdimen\wgap
\newdimen\egap

\def\drawbox#1#2#3{\vbox{
  \setbox\drawing=\vbox{\offinterlineskip\epsfbox{#2.eps}\kern 0pt}
  \drawbp=\epsfurx
  \advance\drawbp by-\epsfllx\relax
  \multiply\drawbp by #1
  \divide\drawbp by 100
  \drawx=\drawbp truebp
  \ifdim\drawx>\hsize\drawx=\hsize\fi
  \epsfxsize=\drawx
  \setbox\drawing=\vbox{\offinterlineskip\epsfbox{#2.eps}\kern 0pt}
  \drawx=\wd\drawing
  \drawy=\ht\drawing
  \ngap=0pt \sgap=0pt \wgap=0pt \egap=0pt 
  \setbox0=\vbox{\offinterlineskip
    \box\drawing \ifgridlines\drawgrid\drawx\drawy\fi #3}
  \kern\ngap\hbox{\kern\wgap\box0\kern\egap}\kern\sgap}}

\def\draw#1#2#3{
  \setbox\drawing=\drawbox{#1}{#2}{#3}
  \advance\figcount by 1
  \goodbreak
  \midinsert
  \centerline{\ifgridlines\boxgrid\drawing\fi\box\drawing}
  \smallskip
  \vbox{\offinterlineskip
    \centerline{Figure~\number\figcount}
    \smash{\marginlabel{#2}}}
  \endinsert}

\def\nextfigtoks{%
  \advance\figcount by 1%
  \edef\numtoks{\number\figcount}%
  \advance\figcount by -1}

\newif\ifgridlines
\newbox\figtbox
\newbox\figgbox
\newdimen\figtx
\newdimen\figty

\newdimen\bwd
\bwd=2sp % 2sp (1/32768") is smallest visible width for Textures

\def\hline#1{\vbox{\smash{\hbox to #1{\leaders\hrule height \bwd\hfil}}}}

\def\vline#1{\hbox to 0pt{%
  \hss\vbox to #1{\leaders\vrule width \bwd\vfil}\hss}}

\def\clap#1{\hbox to 0pt{\hss#1\hss}}
\def\vclap#1{\vbox to 0pt{\offinterlineskip\vss#1\vss}}

\def\hstutter#1#2{\hbox{%
  \setbox0=\hbox{#1}%
  \hbox to #2\wd0{\leaders\box0\hfil}}}

\def\vstutter#1#2{\vbox{
  \setbox0=\vbox{\offinterlineskip #1}
  \dp0=0pt
  \vbox to #2\ht0{\leaders\box0\vfil}}}

\def\crosshairs#1#2{
  \dimen1=.002\drawx
  \dimen2=.002\drawy
  \ifdim\dimen1<\dimen2\dimen3\dimen1\else\dimen3\dimen2\fi
  \setbox1=\vclap{\vline{2\dimen3}}
  \setbox2=\clap{\hline{2\dimen3}}
  \setbox3=\hstutter{\kern\dimen1\box1}{4}
  \setbox4=\vstutter{\kern\dimen2\box2}{4}
  \setbox1=\vclap{\vline{4\dimen3}}
  \setbox2=\clap{\hline{4\dimen3}}
  \setbox5=\clap{\copy1\hstutter{\box3\kern\dimen1\box1}{6}}
  \setbox6=\vclap{\copy2\vstutter{\box4\kern\dimen2\box2}{6}}
  \setbox1=\vbox{\offinterlineskip\box5\box6}
  \smash{\vbox to #2{\hbox to #1{\hss\box1}\vss}}}

\def\boxgrid#1{\rlap{\vbox{\offinterlineskip
  \setbox0=\hline{\wd#1}
  \setbox1=\vline{\ht#1}
  \smash{\vbox to \ht#1{\offinterlineskip\copy0\vfil\box0}}
  \smash{\vbox{\hbox to \wd#1{\copy1\hfil\box1}}}}}}

\def\drawgrid#1#2{\vbox{\offinterlineskip
  \dimen0=\drawx
  \dimen1=\drawy
  \divide\dimen0 by 10
  \divide\dimen1 by 10
  \setbox0=\hline\drawx
  \setbox1=\vline\drawy
  \smash{\vbox{\offinterlineskip
    \copy0\vstutter{\kern\dimen1\box0}{10}}}
  \smash{\hbox{\copy1\hstutter{\kern\dimen0\box1}{10}}}}}

\def\figtext#1#2#3#4#5{
  \setbox\figtbox=\hbox{#5}
  \dp\figtbox=0pt
  \figtx=-#3\wd\figtbox \figty=-#4\ht\figtbox
  \advance\figtx by #1\drawx \advance\figty by #2\drawy
  \dimen0=\figtx \advance\dimen0 by\wd\figtbox \advance\dimen0 by-\drawx
  \ifdim\dimen0>\egap\global\egap=\dimen0\fi
  \dimen0=\figty \advance\dimen0 by\ht\figtbox \advance\dimen0 by-\drawy
  \ifdim\dimen0>\ngap\global\ngap=\dimen0\fi
  \dimen0=-\figtx
  \ifdim\dimen0>\wgap\global\wgap=\dimen0\fi
  \dimen0=-\figty
  \ifdim\dimen0>\sgap\global\sgap=\dimen0\fi
  \smash{\rlap{\vbox{\offinterlineskip
    \hbox{\hbox to \figtx{}\ifgridlines\boxgrid\figtbox\fi\box\figtbox}
    \vbox to \figty{}
    \ifgridlines\crosshairs{#1\drawx}{#2\drawy}\fi
    \kern 0pt}}}}

% macros to add space to text on specified sides

\def\hpad#1#2#3{\hbox{\kern #1\hbox{#3}\kern #2}}
\def\vpad#1#2#3{\setbox0=\hbox{#3}\dp0=0pt\vbox{\kern #1\box0\kern #2}}

% macro to give one text string the apparent height of another

% macro to center one text string over another

\def\stack#1#2#3{\vbox{\offinterlineskip
  \setbox2=\hbox{#2}
  \setbox3=\hbox{#3}
  \dimen0=\ifdim\wd2>\wd3\wd2\else\wd3\fi
  \hbox to \dimen0{\hss\box2\hss}
  \kern #1
  \hbox to \dimen0{\hss\box3\hss}}}

% macros to hide size of trailing exponents

\def\hexp#1{%
  \setbox0=\hbox{${}^{#1}$}%
  \hbox to .5\wd0{\box0\hss}}

%=========================================================================%
% macros for matrices and arrows
%=========================================================================%

% typical usage:
%   \rightarrowmat{2pt}{4pt}{d & bd \cr \!-c & 0 \cr 0 & -ac \cr}

\def\bmatrix#1#2{{\smallmath\left[\vcenter{\halign
  {&\kern#1\hfil$##\mathstrut$\kern#1\cr#2}}\right]}}

\def\rightarrowmat#1#2#3{
  \setbox1=\hbox{\kern#2$\bmatrix{#1}{#3}$\kern#2}
  \,\vbox{\offinterlineskip\hbox to\wd1{\hfil\copy1\hfil}
    \kern 3pt\hbox to\wd1{\rightarrowfill}}\,}

\def\leftarrowmat#1#2#3{
  \setbox1=\hbox{\kern#2$\bmatrix{#1}{#3}$\kern#2}
  \,\vbox{\offinterlineskip\hbox to\wd1{\hfil\copy1\hfil}
    \kern 3pt\hbox to\wd1{\leftarrowfill}}\,}

\def\rightarrowbox#1#2{
  \setbox1=\hbox{\kern#1\hbox{\smallmath #2}\kern#1}
  \,\vbox{\offinterlineskip\hbox to\wd1{\hfil\copy1\hfil}
    \kern 3pt\hbox to\wd1{\rightarrowfill}}\,}

\def\leftarrowbox#1#2{
  \setbox1=\hbox{\kern#1\hbox{\smallmath #2}\kern#1}
  \,\vbox{\offinterlineskip\hbox to\wd1{\hfil\copy1\hfil}
    \kern 3pt\hbox to\wd1{\leftarrowfill}}\,}

%=========================================================================%
% quire macros for preview mode and making booklets
%=========================================================================%

% \legalbooklet{20} makes a booklet from legal paper in landscape
% orientation, where "20" is the page count. To preview, give a negative
% pagecount. Either print using the legal duplex option on a modern laser
% printer, or struggle to simulate this effect manually. Bind using a long
% reach stapler.

% \preview squeezes two pages side by side in landscape orientation. It
% is not suitable for printing, but ideal for previewing on a two page
% monitor.

% \twoup squeezes two pages onto letter paper in landscape mode,
% suitable for printing.

% Each of these macros calls the file "quire.tex"

\def\bookletdims{
  \hsize=5.25truein
  \vsize=7truein
}

\def\legalbooklet#1{
  \input quire
  \bookletdims
  \htotal=7.0truein
  \vtotal=8.5truein
  % below computed from above
  \hoffset=\htotal
  \advance\hoffset by -\hsize
  \divide\hoffset by 2
  \voffset=\vtotal
  \advance\voffset by -\vsize
  \divide\voffset by 2
  \advance\voffset by -.0625truein
  \shhtotal=2\htotal
  % below doesn't need to change
  \horigin=0.0truein
  \vorigin=0.0truein
  \shstaplewidth=0.01pt
  \shstaplelength=0.66truein
  \shthickness=0pt
  \shoutline=0pt
  \shcrop=0pt
  \shvoffset=-1.0truein
  \ifnum#1>0\quire{#1}\else\qtwopages\fi
}

\def\preview{
  \input quire
  \bookletdims
  \hoffset=0.1truein
  \vtotal=8.5truein
  \shhtotal=14truein
  % below computed from above
  \voffset=\vtotal
  \advance\voffset by -\vsize
  \divide\voffset by 2
  \advance\voffset by -.0625truein
  \htotal=2\hoffset
  \advance\htotal by \hsize
  % below doesn't need to change
  \horigin=0.0truein
  \vorigin=0.0truein
  \shstaplewidth=0.5pt
  \shstaplelength=0.5\vtotal
  \shthickness=0pt
  \shoutline=0pt
  \shcrop=0pt
  \shvoffset=-1.0truein
  \qtwopages
}

\def\twoup{
  \input quire
  \hsize=4.79452truein % 5.25/1.095
  \vsize=7truein
  \vtotal=8.5truein
  \shhtotal=11truein
  % below computed from above
  \hoffset=-2\hsize
  \advance\hoffset by \shhtotal
  \divide\hoffset by 6
  \voffset=\vtotal
  \advance\voffset by -\vsize
  \divide\voffset by 2
  \advance\voffset by -12truept
  \htotal=2\hoffset
  \advance\htotal by \hsize
  % below doesn't need to change
  \horigin=0.0truein
  \vorigin=0.0truein
  \shstaplewidth=0.01pt
  \shstaplelength=0pt
  \shthickness=0pt
  \shoutline=0pt
  \shcrop=0pt
  \shvoffset=-1.0truein
  \qtwopages
}

%=========================================================================%
% timestamp (adapted from eplain.tex)
%=========================================================================%

\newcount\countA
\newcount\countB
\newcount\countC

\def\monthname{\begingroup
  \ifcase\number\month
    \or January\or February\or March\or April\or May\or June\or
    July\or August\or September\or October\or November\or December\fi
\endgroup}

\def\dayname{\begingroup
  \countA=\number\day
  \countB=\number\year
  \advance\countA by 0 % adjust after each leap day
  \advance\countA by \ifcase\month\or
    0\or 31\or 59\or 90\or 120\or 151\or
    181\or 212\or 243\or 273\or 304\or 334\fi
  \advance\countB by -1995
  \multiply\countB by 365
  \advance\countA by \countB
  \countB=\countA
  \divide\countB by 7
  \multiply\countB by 7
  \advance\countA by -\countB
  \advance\countA by 1
  \ifcase\countA\or Sunday\or Monday\or Tuesday\or Wednesday\or
    Thursday\or Friday\or Saturday\fi
\endgroup}

\def\timename{\begingroup
   \countA = \time
   \divide\countA by 60
   \countB = \countA
   \countC = \time
   \multiply\countA by 60
   \advance\countC by -\countA
   \ifnum\countC<10\toks1={0}\else\toks1={}\fi
   \ifnum\countB<12 \toks0={\sevenrm AM}
     \else\toks0={\sevenrm PM}\advance\countB by -12\fi
   \relax\ifnum\countB=0\countB=12\fi
   \hbox{\the\countB:\the\toks1 \the\countC \thinspace \the\toks0}
\endgroup}

\def\timestamp{\dayname, \the\day\ \monthname\ \the\year, \timename}

%==========================================================================
% macros (specific to this paper)
%==========================================================================

% surround with $ $ if not already in math mode
\def\enma#1{{\ifmmode#1\else$#1$\fi}}

\def\mathbb#1{{\bbold #1}}
\def\mathbf#1{{\bf #1}}

% blackboard bold symbols

% caligraphic symbols

% bold symbols

% misc

\def\H{\widetilde{H}}

\centerline{\hd  The moduli space of (1,11)-polarized abelian surfaces is 
unirational}
\medskip
\centerline{\it Mark Gross\footnote{*}{Supported by NSF grant DMS-9700761.}
and Sorin Popescu\footnote{**}
{Partially supported by NSF grant DMS-9610205 and MSRI, Berkeley.}}
\medskip
%\centerline{\timestamp}
\medskip
{\settabs 3 \columns
\+Mathematics Institute&&
Department of Mathematics\cr
\+University of Warwick&&
Columbia University\cr
\+Coventry, CV4 7AL, UK&&
New York, NY 10027\cr
\+mgross@maths.warwick.ac.uk&&
psorin@math.columbia.edu\cr}

{\vglue .5in
\narrower{\noindent
\abst Abstract.
We prove that the moduli space $\A_{11}^{lev}$ of 
$(1,11)$-polarized abelian surfaces with level structure
of canonical type is  birational to  Klein's cubic 
hypersurface in $\Pfour$.
Therefore, $\A_{11}^{lev}$ is unirational but not
rational, and there are no $\Gamma_{11}$-cusp forms
of weight 3. The same methods also provide an easy proof of 
the rationality of $\A_{9}^{lev}$.}
\vglue .5in}

%\showlabelsabove        % this prints the labels above each def, thm, etc
%\showlabels            % this prints the labels on the side 
%\preview               % preview format 2 pages side by side

% Forward references 
\forward{torus}{Lemma}{3.4}
\forward{lin.eq}{Theorem}{2.6}

Classical results of Tai, Freitag and Mumford and newer results
of O'Grady, Gritsenko, Hulek and Sankaran say that moduli spaces of 
polarized abelian varieties are almost always of general type. 
However, for abelian varieties of small dimension and  
polarizations of small degree 
the situation is different and the corresponding moduli spaces 
usually have beautiful geometry.

In this paper we describe a projective model for the moduli
of complex abelian surfaces with a polarization of type $(1,11)$, with
level structure of canonical type. As a direct consequence we obtain the 
unirationality of this moduli space, which also turns out to
be non-rational.  However, unirationality already implies that 
there exist no $\Gamma_{11}$-cusp forms
of weight 3.

Let $\A_d$ denote the moduli space of polarized abelian surfaces
of type $(1,d)$, and let $\A_d^{lev}$ be the moduli space of
$(1,d)$-polarized abelian surfaces with canonical level structure. 
The map which forgets the level structure
represents $\A_d^{lev}$ as a finite cover of  $\A_d$. Its general
fiber is $\Gamma_d/\Gamma_d^{lev}\cong \SL_2(\boldz_{d})$, where
$\Gamma_d$ and $\Gamma_d^{lev}$ denote the corresponding paramodular
groups. In particular, if $d$ is an odd prime number, then the 
forgetful morphism is a ramified cover of degree $d(d^2-1)/2$.
(See [LB], [Mum], [GP1] for the definition of canonical level structure, 
and basic results.)

Our main result is the following 

\theorem{main} The moduli space $\A_{11}^{lev}$ is birational to  
Klein's cubic hypersurface
$$\K=V(\sum_{i\in\boldz_5} x_i^2x_{i+1}=0)\subset\Pfour.$$ 
In particular, $\A_{11}^{lev}$ is unirational but not
rational. 

The cubic hypersurface $\K\subset\P^4$ was first studied by Klein [Kl] 
(see also [KlF], Band II) in connection with the $z$-embedding of 
the modular curve $X(11)$ of level 11, 
which turns out to be defined by the $4\times 4$-minors of 
the Hessian of the equation of $\K$. In
this respect, we note that $\K$ is the unique 
$\PSL_2(\boldz_{11})$-invariant of degree three
in $\Pfour$, and that furthermore 
$\PSL_2(\boldz_{11})$ is its full automorphism group [Ad1]. 
The Klein cubic being smooth is unirational but not rational, cf. 
[CG], [Mur], [Bea].

Our result should be regarded in light of the following facts:
${\A}_t$ is not unirational (and in fact $p_g(\widetilde{\A}_t)\ge 1$) 
if $t\ge 13$ and 
$t\ne 14,\; 15,\; 16,\; 18,\; 20,\; 24,\; 30,$ $36$ 
(Gritsenko [Gri1], [Gri2]), while
$\widetilde\A^{lev}_p$ is a 3-fold of general type for all primes numbers
$p\ge 37$ (Hulek-Sankaran [HS1], Gritsenko-Hulek, appendix to [Gri1]),
where $\widetilde\A^{lev}_p$ is a smooth projective model of a 
compactification of $\A^{lev}_p$.
See also [GH], and the survey paper [HS2] for related results, and
[Bo] for a finiteness result in the same spirit. On the other hand, 
$\overline{\A^{lev}_5}\cong\P(H^0(F_{HM}(3)))$, 
where bar stands for the Igusa (=Voronoi) toroidal compactification and 
$F_{HM}$ is the Horrocks-Mumford bundle on $\P^4$ ([HM], [HKW]), 
while $\A^{lev}_7$ is rational having as birational model a 
smooth $V_{22}$, a  prime Fano 3-fold of index 1 and  genus 12 
which is rational (see [MS], [Schr] and [GP2] for details). 
It would be interesting
to know how the Klein cubic ``compares'' with the toroidal compactification
of $\A^{lev}_{11}$. 

In a series of forthcoming papers [GP2], [GP3], we will give details 
as to the structure of $\A_d^{lev}$, $6\le d\le 12$, (excluding
$d=9$ and $11$, which are covered here) and $\A_d$, 
$d=14, 16, 18$ and $20$. In particular, we will prove 
their rationality or unirationality.

Finally, the methods used in this paper also provide an easy proof of 
the rationality (over $\QQ(\xi)$ for $\xi$ a primitive $9^{\rm th}$ 
root of unity) of $\A^{lev}_9$. The unirationality of
this space also follows implicitly from O'Grady's work [O'G]. 
He identifies
$\A^{lev}_{p^2}$, for $p$ prime, 
with the moduli space $\A_1(p)$ of pairs of principally polarized abelian 
surfaces and rank two subspaces of the $p$-torsion points, 
non-isotropic for the Weil pairing.
O'Grady studies the extension to (natural) toroidal compactifications
of the finite natural forgetful map $\pi$ from $\A_1(p)$ to the 
moduli space $\A_1$ of
principally polarized abelian surfaces.  Not all singularities of 
the toroidal compactification
of $\A_1(p)$ are canonical, so O'Grady needs to describe carefully
a partial desingularization all of whose singularities are canonical, before
being able to apply Hurwitz's formula for $\pi$ to get an expression for
the canonical class.   For $p=3$, our method in addition to being simpler
also has the advantage of providing an explicit
rational parametrization (over $\QQ(\xi)$).

{\it Acknowledgments:} We thank Igor Dolgachev, David Eisenbud, 
Klaus Hulek, and Kristian Ranestad for many useful
discussions, and Allan Adler and Gregory Sankaran for a careful
reading of a preliminary version of this paper. 
We are also grateful to Dave Bayer, Dan Grayson and 
Mike Stillman for {\sl Macaulay} [BS], and {\sl Macaulay2} [GS] 
which helped us tremendously to understand the
shape of the equations described in this paper. The second author also thanks
the Mathematical Sciences Research Institute, Berkeley for its hospitality
while part of this paper was being written.

\section {representations} {Preliminaries}

We review basic properties of polarized abelian surfaces, specializing
to the case of polarizations of type $(1,11)$. For a more detailed
review of this material as used in this paper, see [GP1], \S 1. In general, we
use the notation and definitions of [LB] and [Mum].

Let $(A,\L)$ be a general abelian surface with a polarization of type $(1,11)$.
Then $|\L|$ induces an embedding of $A\subset\P^{10}=\P(\dual{H^0(\L)})$
as a projectively normal surface  
(cf. [Laz]) of degree 22 and sectional genus 12.
(The projective normality of the general such abelian surface, 
follows also from the proof of [GP1], Theorem 6.5.)  
Riemann-Roch tells us that $A$ is contained in $22$ quadrics, 
which generate the homogeneous ideal of $A$ ([GP1], Theorem 6.5).

The line bundle $\L$ induces a natural
map from $A$ to its dual, $\phi_{\L}:A\rightarrow \hat A$,
given by $x\mapsto t_x^*\L\otimes\L^{-1}$,
where $t_x:A\rightarrow A$ is the morphism given by translation by $x\in A$. 
Its kernel $K(\L)$ is isomorphic to $\boldz_{11}\times\boldz_{11}$,
and is dependent only on the polarization. 

For every $x\in K(\L)$ there is an isomorphism $t_x^*\L
\cong \L$. This induces a projective representation 
$K(\L)\rightarrow {\rm PGL}(H^0(\L))$, which lifts uniquely to a 
linear representation of $K(\L)$ after taking a central extension 
of $K(\L)$ 
$$1 \rTo\, {{\bf C}^*} \rTo\, {\G(\L)} \rTo\, {K(\L)} \rTo\, 0,$$
whose Schur commutator map is the Weil pairing. 
$\G(\L)$ is the {\it theta group} of $\L$ and is isomorphic to the abstract 
Heisenberg group $\H(11)$, while the above representation 
is isomorphic to the Schr\"odinger representation of $\H(11)$ on
$V={\bf C}(\boldz_{11})$,
the vector space of complex-valued functions on  
$\boldz_{11}$.
An isomorphism between $\G(\L)$ and $\H(11)$, which restricts 
to the identity on centers induces a symplectic isomorphism 
between $K(\L)$ and $\boldz_{11}\times\boldz_{11}$. Such an
isomorphism is called 
a {\it level structure of canonical type} on $(A,c_1(\L))$. (See [LB],
Chapter 8, \S 3 or [GP1], \S 1.)

A decomposition $K(\L)=K_1(\L)\oplus K_2(\L)$,
with $K_1(\L)\cong K_2(\L) \cong \boldz_{11}$ subgroups 
isotropic with respect to the Weil pairing, and a choice
of a characteristic $c$ ([LB], Chapter 3, \S 1) for $\L$, define a unique basis
$\{\vt^c_x \mid x\in K_1(\L)\}$ of {\it canonical theta functions} 
for the space $H^0(\L)$ defined in [LB], Chapter 3, \S 2. This
basis allows an identification of $H^0(\L)$ with $V$ via 
$\vt_{\gamma}^c\mapsto x_{\gamma}$, where $x_{\gamma}$
is the function on $\boldz_{11}$ defined by 
$x_{\gamma}(\delta)=\cases{1&$\gamma=\delta$\cr 0&$\gamma\not=\delta$\cr}$
for $\gamma,\delta\in\boldz_{11}$. The $x_0,\ldots,x_{10}$ can also be
identified with coordinates on $\P(\dual{H^0(\L)})$. 
Under this identification, the
representation $\G(\L)\rightarrow {\rm GL}(H^0(\L))$ 
coincides with the  Schr\"odinger representation 
$\H(11)\rightarrow GL(V)$. We will only consider the
action of $\HHH_{11}$,  the finite subgroup of $\H(11)\rightarrow GL(V)$ 
generated in the Schr\"odinger representation by $\sigma$ and $\tau$, where 
$$ \sigma(x_i)=x_{i-1}, \qquad \tau(x_i)={\xi}^{-i} x_i,$$
for all $i\in{\boldz}_{11}$ and $\xi=e^{2{\pi}i\over 11}$ is a primitive
root of unity of order 11. Notice that 
$\lbrack\sigma,\tau\rbrack=\xi$,
so ${\HHH}_{11}$ is a central extension 
$$1 \rTo\; {\bf\mu_{11}} \rTo\;  {\HHH_{11}} \rTo\; 
{\boldz_{11}}\times{\boldz_{11}} \rTo\; 0.$$
Thus the choice of a canonical level structure means that 
if $A$ is embedded in $\P(\dual{H^0(\L)})$ using as coordinates
$x_{\gamma}=\vt_{\gamma}^c$,  $\gamma\in
\boldz_{11}$, then the image of $A$ will be invariant under the action of the
Heisenberg group $\HHH_{11}$ via the Schr\"odinger representation. (See
[LB], Chapter 6, \S 7).
 
If moreover the line bundle $\L$ is chosen to be symmetric 
(and there are always finitely many choices of such
an $\L$ for a given polarization type), then the  embedding via $|\L|$ 
is also invariant under the involution $\iota$, where
$$ \iota(x_i)=x_{-i}, \qquad i\in{\boldz}_{11},$$
which restricts to $A$ as the involution $x\mapsto -x$. 

Let $N(\HHH_{11})$ be the normalizer of $\HHH_{11}$ inside $\SL(V)$, 
where the inclusion $\HHH_{11}\subset \SL(V)$ is 
the Schr\"odinger representation. An element $\alpha\in N(\HHH_{11})$ 
induces an outer automorphism of $\HHH_{11}$, and hence an 
automorphism of $\boldz_{11}\times\boldz_{11}$
preserving the Weil pairing $e^D$, for $D=(1,11)$.
The group of such automorphisms is 
$\SL_2(\boldz_{11})$, and thus
we get a map $\psi : N(\HHH_{11})\rightarrow \SL_2(\boldz_{11})$. 
As in [HM] \S 1, one sees that the kernel of this map is $\HHH_{11}$ and that 
$\psi$ is surjective. This leads to extensions
\newarrow{Equals}=====
$$\diagram[midshaft,small]
1&\rTo &\HHH_{11}&\rTo &N(\HHH_{11})&\rTo^\psi& \SL_2(\boldz_{11})&\rTo& 1\\
 &     &\dTo     &     &\dTo        &    & \dEquals \\
0&\rTo &\boldz_{11}\times\boldz_{11} &\rTo &N(\HHH_{11})/{\bf{\mu_{11}}}
&\rTo& \SL_2(\boldz_{11})&\rTo& 1\\
\enddiagram$$
where in the bottom row  $N(\HHH_{11})/{\bf{\mu_{11}}}$ is
a semidirect product by the above (symplectic) 
action of $\SL_2(\boldz_{11})$. Since $H^2(\SL_2(\boldz_{11}),{\bf C}^*)=0$
it follows  that $N(\HHH_{11})$ is in fact the semi-direct product
$\HHH_{11}\rtimes\SL_2(\boldz_{11})$ (see [HM], \S 1 for details in
the identical case of $\HHH_5$). 

Therefore the Schr\"odinger representation of $\HHH_{11}$ 
induces an 11 dimensional representation
$$\rho_{11}:\SL_2(\boldz_{11})\rTo \SL(V).$$
In terms of generators and relations, cf. [BM], one has:
$${\rm PSL}_2(\boldz_{11})=\langle S,T\mid S^{11}={(S^2TS^6T)}^3=1,\quad {(ST)}^3=T^2=1\rangle,$$
where 
$$S=\pmatrix{1&1\cr 0&1\cr}\qquad
{\rm and}\qquad T=\pmatrix{0&\hfill -1\cr 1&\hfill 0\cr},$$
while the representation $\rho_{11}$ is given projectively by
$$\rho_{11}(S)=(\xi^{ij/2}\delta_{ij})_{0\le i,j \le 10}\qquad
{\rm and}\qquad
\rho_{11}(T)={1\over\sqrt{11}}(\xi^{-ij})_{0\le i,j \le 10},$$
where $\xi$ is the above fixed $11^{\rm th}$ root of unity. 
(See [Tan] and [Si] for details.)

The center of $\SL_2(\boldz_{11})$ is generated by $T^2$, 
and $\rho_{11}(T^2)=-\iota$. 
Thus the representation $\rho_{11}$ is reducible. In fact, if $V_+$ and $V_-$
are the positive and negative eigenspaces, respectively,  
of the involution $\iota$ acting on $V$, then $V_+$ and $V_-$ 
are easily seen to be invariant under $\rho_{11}$,
and moreover $\rho_{11}$ splits as $\rho_+\oplus \rho_-$, where
$\rho_{\pm}$ is the representation of $\SL_2(\boldz_{11})$ acting on
$V_{\pm}$. Note that $\rho_-$ is trivial
on the center of $\SL_2(\boldz_{11})$, so it descends to give an irreducible 
representation
$$\rho_-:\PSL_2(\boldz_{11})\rTo \GL(V_-).$$
For the reader's convenience, we reproduce from the Atlas 
of finite groups [CNPW] the character table for $\PSL_2(\boldz_{11})$:
$$
\vbox{\tabskip=0pt \offinterlineskip
\halign to 400pt
{\hfil#$\;$\vrule&\strut\hfil$\quad\quad$#\hfil&\hfil$\quad\quad$#\hfil&\hfil$\quad\quad$#
\hfil&\hfil$\quad\quad$#\hfil&\hfil$\quad\quad$#\hfil&
\hfil$\quad\quad$#\hfil&\hfil$\quad\quad$#\hfil&\hfil$\quad\quad$#\hfil\cr
Size of conjugacy class&1&55&110&132&132&110&60&60\cr
Conjugacy class&$I$&
$\gamma_1$&$\gamma_2$&$\gamma_3$&$\gamma_4$&
$\gamma_5$&$\gamma_6$&$\gamma_7$\cr
Character&&&&&&&&\cr
\noalign{\hrule}
$\chi_1$&1&1&1&1&1&1&1&1\cr
$\chi_2$&5&1&-1&0&0&1&$\beta$&$\bar\beta$\cr
$\chi_3$&5&1&-1&0&0&1&$\bar\beta$&$\beta$\cr
$\chi_4$&10&-2&1&0&0&1&-1&-1\cr
$\chi_5$&10&2&1&0&0&-1&-1&-1\cr
$\chi_6$&11&-1&-1&1&1&-1&0&0\cr
$\chi_7$&12&0&0&$\alpha$&$\alpha'$&0&1&1\cr
$\chi_8$&12&0&0&$\alpha'$&$\alpha$&0&1&1\cr}
}
$$
where $\alpha={1\over 2}(-1+\sqrt{5}), \alpha'={1\over 2}(-1-\sqrt{5}),$
and $\beta={1\over2}(-1+\sqrt{-11})$. The conjugacy classes are 
represented by
$$
\vbox{\tabskip=0pt \offinterlineskip
\halign to 400pt
{\strut\hfil#\quad \hfil&\hfil#\quad \hfil&\hfil#\quad
\hfil&\hfil#\quad \hfil&\hfil#\quad \hfil&
\hfil#\quad \hfil&\hfil#\quad \hfil&\hfil#\hfil\cr
$I$&$\gamma_1$&$\gamma_2$&$\gamma_3$&$\gamma_4$&
$\gamma_5$&$\gamma_6$&$\gamma_7$\cr
$\pmatrix{1&0\cr 0&1}$&$\pmatrix{0&-1\cr 1& 0\cr}$&$\pmatrix{
1& -1\cr 1& 0\cr}$&$\pmatrix{3&0\cr 0&4\cr}$
&$\pmatrix{5&0\cr 0&9\cr}$&$\pmatrix{3&2\cr 4&3\cr}$&$\pmatrix{1&1\cr
0&1\cr}$&$\pmatrix{1&2\cr 0&1\cr}$\cr}}$$

Similarly, the other direct summand 
$$\rho_+:\SL_2(\boldz_{11})\rTo \GL(V_+)$$
is one of two mutually dual $6$-dimensional 
irreducible representations of $\SL_2(\boldz_{11})$. 
We refer the reader to [Dor] for the character table of $\SL_2(\boldz_{11})$.

\section {moduli11} {Moduli of $(1,11)$-polarized abelian surfaces}

From now on, let $(A,\L)$ be a general abelian surface with a 
polarization of type $(1,11)$
and with canonical level structure. As seen in the previous section,
$|\L|$ embeds $A\subset\P^{10}=\P(\dual{H^0(\L)})=\P(\dual{V})$  
as a projectively normal surface of degree 22 and sectional genus 12,
which is invariant under the action of the
Heisenberg group $\HHH_{11}$ via the Schr\"odinger representation. 
In particular, $H^0(\I_A(n))$ is also a 
representation of weight $n$ of the Heisenberg group (i.e., a central element 
$z\in{\bf C}^{\ast}$ acts by multiplication with $z^n$), 
and hence all its irreducible components will
have dimension $11/\gcd(11,n)$. (See [LB], pg. 179, for this last fact.)

We will first determine equations for the locus 
of odd two-torsion points of $(1,11)$ polarized abelian surfaces. This
is a set 
in $\P^-=\P(\dual{V_-})$. 
To analyze the equations which arise, 
we will need to make use of the $\SL_2(\boldz_{11})$ symmetry present.

Notice that $S^2(V)=H^0(\O_{\P(\dual{H^0(\L)})}(2))$ is $66$-dimensional, and
as a representation of weight $2$,  splits
into six isomorphic 11 dimensional representations of
$\HHH_{11}$, each isomorphic to a twist $V'$ of
the Schr\"odinger representation. On the other hand, by (2.2) of [HM], 
$\hom_{\HHH_{11}}(V',V\otimes V)$ is an 
$N(\HHH_{11})/\HHH_{11}\cong \SL_2(\boldz_{11})$-module
and coincides with the representation $\rho_{11}$ described
in \ref{representations} Decomposing into the positive and 
negative eigenspaces of the involution $\iota$
we deduce that
$$V_+=\hom_{\HHH_{11}}(V',S^2(V))\quad {\rm and}\quad
V_-=\hom_{\HHH_{11}}(V',\wedge^2(V)),$$ 
as $\SL_2(\boldz_{11})$-modules. By (2.2) of [HM],
$\varphi: V'\otimes \hom_{\HHH_{11}}(V', S^2(V))\rightarrow S^2(V)$ 
is an isomorphism of $N(\HHH_{11})$-modules, and since
$S^2(V)$ is an irreducible $N(\HHH_{11})$-module it follows
that $V_+$ is a six dimensional irreducible representation 
of $\SL_2(\boldz_{11})$.

The 
$\SL_2(\boldz_{11})$-isomorphism  $\varphi: V_+\otimes V'\rightarrow S^2(V)$
can be represented as follows. We use the usual basis $x_0,\ldots,x_{10}$
for $V$ (identified with the basis of canonical theta functions of $H^0(\L)$),
and the basis $f_0,\ldots,f_{10}$ of $V'$ such that $\sigma(f_i)=f_{i-1}$
and $\tau(f_i)=\xi^{-2i}f_i$, $i\in \boldz_{11}$. 
Then there is a basis $e_0,\ldots,e_5$
of $V_+$ (in fact the projection of $x_0,\ldots,x_5$ onto $V_+$
coming from the decomposition $V=V_+\oplus V_-$) such that the map
$\varphi$ takes $e_i\otimes f_j$ to the $(i,j)^{\rm th}$ entry of the
$6\times 11$ matrix $R_5$, whose entries are 
$${(R_5)}_{ij}=x_{j+i}x_{j-i},\quad 0\le i \le 5, 0\le j \le 10,$$
where the indices of the variables are $\mod 11$. Thus the span of the entries
in any row of $R_5$ are a $\HHH_{11}-$subrepresentation of $S^2(V)$.
Also, any $\HHH_{11}-$subrepresentation of $S^2(V)$ can
be obtained by taking a linear combination of the rows, and taking the
span in $S^2(V)$ of the resulting 11 quadratic polynomials. In addition,
if $P\in\P^{10}$ and $v\in V_+$, then $v\cdot R_5(P)=0$ if and only
if $P$ is contained in the scheme cut out by the $\HHH_{11}$-subrepresentation
of quadrics determined by $v$.

Abusing
notation, we will also denote by $V_\pm$ the eigenspaces 
of the involution $\iota$ acting on $V'$ (by $\iota(f_i)=f_{-i}$). 
Then the restriction of $\varphi$ to $V_+\otimes V_+$
induces a $\SL_2(\boldz_{11})$-isomorphism 
$$\Phi : \wedge^2(V_+)\rTo S^2(V_-),$$
usually called the {\it intertwining operator} (see [We],
[AR] pp. 62--63, 74, or [Ad5] for details). We may regard
the intertwining isomorphism  $\Phi$ as being induced by a skew-symmetric
matrix with entries quadratic polynomials in the coordinates
of $V_-$, namely: $\Phi$ takes the element 
$e_i\wedge e_j$ of $\wedge^2(V_+)$
to the $(i,j)^{\rm th}$ entry of the matrix
$$S=
\pmatrix{0&{{{x}}_{1}}^{{2}}&{{{x}}_{{2}}}^{{2}}&{{{x}}_{{3}}}^{{2}}&
{{{x}}_{{4}}}^{{2}}&{{{x}}_{{5}}}^{{2}}\cr 
{-{{{x}}_{1}}^{{2}}}&0&{{x}}_{1}{{x}}_{{3}}&{{x}}_{{2}}{{x}}_{{4}}&
{{x}}_{{3}}{{x}}_{{5}}&{-{{x}}_{{4}}{{x}}_{{5}}}\cr 
{-{{{x}}_{{2}}}^{{2}}}&{-{{x}}_{1}{{x}}_{{3}}}&0&{{x}}_{1}{{x}}_{{5}}
&{-{{x}}_{{2}}{{x}}_{{5}}}&{-{{x}}_{{3}}{{x}}_{{4}}}\cr 
{-{{{x}}_{{3}}}^{{2}}}&{-{{x}}_{{2}}{{x}}_{{4}}}&{-{{x}}_{1}{{x}}_{{5}}}
&0&{-{{x}}_{1}{{x}}_{{4}}}&{-{{x}}_{{2}}{{x}}_{{3}}}\cr 
{-{{{x}}_{{4}}}^{{2}}}&{-{{x}}_{{3}}{{x}}_{{5}}}&{{x}}_{{2}}{{x}}_{{5}}
&{{x}}_{1}{{x}}_{{4}}&0&{-{{x}}_{1}{{x}}_{{2}}}\cr 
{-{{{x}}_{{5}}}^{{2}}}&{{x}}_{{4}}{{x}}_{{5}}&{{x}}_{{3}}{{x}}_{{4}}
&{{x}}_{{2}}{{x}}_{{3}}&{{x}}_{1}{{x}}_{{2}}&0\cr }.
$$
Here we are abusing notation by identifying $x_i\in V$ with the
projection of $x_i$ onto $V_-$ under the decomposition $V=V_+\oplus
V_-$.
Note that $S$ can also be viewed as the
restriction of the first $6\times 6$-block of $R_5$  to $\P^-$, where
we use $x_1,\ldots,x_5$ as coordinates on $\P^-$ (where $x_{i}=x_{-i}$). 

We will only need the fact that $S$ arises from the intertwining
operator in Lemma 2.1 1) below. The key fact we will need later
about $S$ is the simple observation that 
if $P\in\P^-$, then $v\cdot R_5(P)=0$ 
if and only if $v\cdot S(P)=0$.

 Following [GP1], \S 6, we define  $D_i\subseteq\P^-$ to
be the locus where the matrix $S$ has rank $\le 2i$, for $i=1,2$.
By the previous observation, we may interpret $D_i$ as the locus
of points in $\P^-\subset\P^{10}=\P(\dual{V})$ 
which are contained in at least
a $(6-2i)$-dimensional family of $\HHH_{11}$-representations of quadrics.
Remark also that all the loci $D_i$ are invariant under the action
of $\PSL_2(\boldz_{11})$ via the representation $\rho_{-}$
defined in \ref{representations}

\lemma{modular} 
\item{1)} $D_1\subset\P^-$ is a smooth curve of degree 20
and genus 26 isomorphic to the modular curve $X(11)$.
(This is Klein's $z$-model of the modular curve $X(11)$. It is the
``trace'' of the origins in the Shioda compactification 
of elliptic normal curves with level structure
in $\P^{10}$, and its embedding is induced by $\lambda^4$, where $\lambda$,
a $10^{\rm th}$ root of the canonical bundle on $X(11)$, is the
generator of the group of $\PSL_2(\boldz_{11})$-invariant 
line bundles on $X(11)$. See
[Kl], pp. 153-156, [AR] and [Dol] for details).
\item{2)} $D_2\subset\P^-$ is an irreducible sextic hypersurface,
defined by the $6\times 6$ Pfaffian of $S$.
It contains as an open subset the locus in 
$\P^{10}$ of odd $2$-torsion points of $(1,11)$-polarized abelian 
surfaces with canonical level structure.

\proof: 
$1)$  Consider the composition
$$\Psi:\P^-=\P(\dual{V_-})\rTo \P(S^2(\dual{V_-}))
\rTo^{\cong} \P(\wedge^2(\dual{V_+})),$$
where the first map is the Veronese embedding and the second
isomorphism is induced by the intertwining operator. Then 
for $P\in\P^-$, one can write $\Psi(P)=S(P)$, the latter being
a skew-symmetric matrix which should be interpreted as an
element of $\wedge^2(\dual{V_+})$. It is then clear that $D_1$
is the pull-back under $\Psi$ of the locus in $\P(\wedge^2(\dual{V_+}))$
of rank 2 skew-symmetric matrices. This latter locus can be identified
with $\Gr(2,\dual{V_+})$. Thus $D_2$ is isomorphic to the pull-back of
$\Gr(2,\dual{V_+})$ under $\Psi$.

Adler-Ramanan [AR], Theorem 19.17, study this pull-back. In particular,
they show that this pull-back gives the same variety in $\P^-$ as
that defined by Klein's equations studied by V\'elu in [Ve].
It is proven in [Ve], (summarized in Th\'eor\`eme  10.6), 
that Klein's equations
give the so-called the $z$-model of $X(11)$, and this model
is always nonsingular. (See [Dol], Theorem 5.1 for the explicit
statement.) Thus $D_1$  is isomorphic
to $X(11)$. 

The genus of $X(11)$ is well-known: see for example [Hu], pg. 59,
for the genus of the modular curve of level $n$. The degree of the
$z$-model is calculated in [Ve]. (See also [AR], Corollary 23.28.)

$2)$ For a general point $P\in\P^-=\P(V_-)$, $S$ has rank 6,
the intertwining operator being an isomorphism, thus
$\P^-=D_3\not=D_2$. In fact $D_2\subset\P^-$ is the sextic hypersurface
given by the Pfaffian of $S$. For the record, the equation $f_6$ of $D_2$  is
$$\eqalign{f_6=-&x_1^2x_2x_3^3+x_1^3x_3x_4^2-x_2^3x_3^2x_5+x_1x_4^3x_5^2
+x_2^2x_4x_5^3\cr
&+x_1x_2^4x_4-x_2x_3x_4^4-x_1^4x_2x_5+x_3^4x_4x_5+x_1x_3x_5^4\cr
&+x_1x_2x_3^2x_4^2-x_1^2x_2^2x_3x_5-x_1x_2^2x_4^2x_5
-x_1^2x_3x_4x_5^2+x_2x_3^2x_4x_5^2,\cr}$$
though we will not make use of its explicit form.
To prove its irreducibility we will use the fact that it
is a $\PSL_2(\boldz_{11})$-invariant in $S^6(V_-)$.
Without loss of generality we may assume that
as a $\PSL_2(\boldz_{11})$ representation 
$V_-$ has character $\chi_3$.
It is easy to see that 
$\chi_{S^2(V_-)}=\chi_3+\chi_5,$
so in particular there are no $\PSL_2(\boldz_{11})$ invariants in $S^2(V_-)$. 
On the other hand, $\chi_{S^3(V_-)}=\chi_1+\chi_5+\chi_6+\chi_7,$
so in $S^3(V_-)$ there is precisely one $\PSL_2(\boldz_{11})$-invariant, which
we will denote by $f_3$. Thus the only way that $f_6$ could fail 
to be irreducible
is if $f_6=f_3^2$. But $f_6$ is not a square. To see this, set, say $x_4=x_5=0$
in the matrix $S$ and take its Pfaffian, or just set $x_4=x_5=0$ in 
the above equation for $f_6$. We get only one term, $-x_1^2x_2x_3^3$, 
which is not a square, so $f_6$ itself cannot be a square, and thus 
$f_6$ is irreducible. (Alternatively, the irreducibility of $D_2$
follows from [Ad3].)

For any point $P\in D_2\setminus D_1$,
$S(P)$ is of rank 4, thus $P$ is contained in a pencil
of $\HHH_{11}$-subrepresentations of quadrics in $\P^{10}$. 
The last claim follows now from the fact  that  
$\HHH_{11}\rtimes\langle\iota\rangle$-invariant abelian surfaces in
$\P^{10}$ lie on a pencil
of $\HHH_{11}$-representations of quadrics, and that odd 2-torsion
points of a general abelian surface get mapped 
to $D_2\setminus D_1\subset\P^-$,
and uniquely determine the surface (cf. [GP1], Lemma 6.3 and Lemma 6.4).\Box

By \ref{modular}, we may define now the morphism 
$$\Theta : D_2\setminus D_1\rTo \Gr(2, V_+)=\Gr(2,6)$$  
$$D_2\setminus D_1\ni P \rMapsto \ker(S(P))=\{v\in V_+ \mid v\cdot S(P)=0\},$$ 
which sends a point $P$ to the pencil of $\HHH_{11}$-subrepresentations
of quadrics containing it.  The general $(1,11)$-polarized abelian
surface embedded with level structure via a symmetric line bundle
meets $\P^-$ in the (images of the) 6 odd 2-torsion points. 
By [GP1], Lemma 6.4 they are mapped via $\Theta$ to a single point in
$\Gr(2,6)$. Thus $\Theta$ factorizes as a rational map               
$$\Theta_{11}:\A_{11}^{lev}\rDashto \Gr(2,6)=\Gr(2, V_+),$$ 
which essentially takes an abelian surface $A$ to
the point in $\Gr(2,6)$ corresponding to the
$\HHH_{11}$-subrepresentation $H^0(\I_A(2))\subset H^0(\O_{\P(V)}(2))$.

\theorem{rat.map11} The map $\Theta_{11}:\A_{11}^{lev}\rDashto
\Gr(2,6)$ yields a birational map between $\A_{11}^{lev}$ and
$\im(\Theta)$.

\proof: This is an immediate consequence of the degeneration arguments in
[GP1], Theorem 6.5,  and the above results. \Box

In order to prove \ref{main} we will need to determine 
the precise structure of 
$\im(\Theta)$. To do so, we will use the following representation of the 
Pl\"ucker embedding of the Grassmannian $\Gr(2,2m)=\Gr(2,W)$, where
$W\cong {\bf C}^{2m}$, $m\ge 2$. 

$\Gr(2,2m)$ is embedded in $\P^{{2m \choose 2}-1}=\P(\wedge^2(W))$ 
via the Pl\"ucker embedding, 
as the variety of those $2$-vectors which are totally decomposable. 
Thus hyperplanes in the Pl\"ucker embedding can be identified 
with (projectivized) skew symmetric
forms $H\in\P(\wedge^2(V)^*)$, and thus with  
$2m\times 2m$ skew-symmetric matrices. In this
setting the Grassmannian $\Gr(2,2m)$ can be also identified 
(as an embedded variety) with the subvariety $R_1$ of 
$2m\times 2m$ skew-symmetric matrices of rank two.
The hyperplane sections corresponding to points of $R_1$ are the 
Schubert cycles $\sigma_1$, the sets of lines intersecting a given 
subspace of codimension 2 in $\P(W)$. 

More precisely, if $P\in\Gr(2,2m)$ corresponds to a subspace of 
$W\cong {\bf C}^{2m}$ spanned by the rows of a $2\times 2m$ matrix
$$L=\pmatrix{a_1&\cdots&a_{2m}\cr
b_1&\cdots&b_{2m}\cr},$$
then the corresponding  $2m\times 2m$ skew-symmetric matrix $H_P=(p_{ij})$ 
has as entries  the Pl\"ucker coordinates of $L$
$$p_{ij}=a_ib_j-a_jb_i.$$
This matrix is rank 2, all the rows being linear combinations
of the rows of $L$. Conversely, given any  $2m\times 2m$ skew-symmetric matrix
of rank 2, the span of the rows yields a two-dimensional subspace of
${\bf C}^{2m}=W$, and hence a point in $\Gr(2,2m)$. Furthermore,
with obvious abuse of notation, the following correspondence holds: 

\lemma{chords} For $H\in\P(\wedge^2(W)^*)$ and $k\in\{1,\ldots,m\}$ the 
following are equivalent
\item{1)} $H\in R_k:=\{H\in \P(\wedge^2(W)^*)\mid \rank(H)\le 2k\}$
\item{2)} $H$ lies in the $k$-chordal locus of $R_1\cong\Gr(2,2m)$ (i.e.,
lies in a linear subspace $\P^{k-1}\subset\P(\wedge^2(W)^*)$ which
is $k$-secant to $R_1$).

\proof: All these facts are classical, and easy to prove.
See for instance [SR].\Box 

An easy computation shows that the codimension of $R_{m-k}$ in
$\P(\wedge^2(W)^*)$ is ${2k\choose 2}$. In particular, 
$R_{m-1}$ is a hypersurface in $\P(\wedge^2(W)^*)$
of degree $m$, defined by the Pfaffian of the generic
skew-symmetric matrix. 

\remark{orbits} The Pl\"ucker embedding is compatible with the 
natural action of $\PGL(W)$, and the orbits under this
action are exactly $R_k\setminus R_{k-1}$, for $k=\overline{1,m}$,
where $R_0=\emptyset$.

The following lemma describes the map $\Theta$ in the above setting:

\lemma{pfaff} Let $M$ be a $2m\times 2m$ skew-symmetric
matrix of forms of rank $2m-2$ on a variety $X$. Then the map
$\Theta:X\rightarrow\Gr(2,2m)$ induced by $x\in X\mapsto 
\ker M(x)\subseteq {\bf C}^{2m}$ is given in (dual) 
Pl\"ucker coordinates by 
$$x\rMapsto M^*(x)\in R_1=\Gr(2,2m),$$
where $M^*$ is the $2m\times 2m$ skew-symmetric matrix defined by
$$M_{ij}^*=\cases{
(-1)^{i+j}\Pf^{ij}(M)& $i<j$,\cr
0&$i=j$,\cr
(-1)^{i+j+1}\Pf^{ij}(M)& $i>j$,\cr}$$
and where $\Pf^{ij}(M)$ is the Pfaffian of the matrix obtained
by deleting the $i^{\rm th}$ and $j^{\rm th}$ 
rows and columns from $M$.

\proof: 
We'll make use of the following standard facts concerning Pfaffian
identities (see [BE] and [Re] for more details). Let $F$ be a free 
module of rank $n=2m$ over the ring $R$, let $F^\ast$ denote the dual module,
and let $f:F^\ast\to F$ be a skew-symmetric morphism 
(that is the  matrix of $f$ corresponding to the choice of a basis
in $F$ and the dual basis in $F^\ast$ is a skew-symmetric
$n\times n$-matrix $M$). Now, $\wedge F^\ast$ and 
$\wedge F$ are modules over each other and we adopt here [BE]'s notation
in writing $a(b)$ for the result of an operation of $a\in\wedge F$
on $b\in\wedge F^\ast$ and vice-versa; thus $a(b)\in\wedge F^\ast$ and 
$b(a)\in\wedge F$. The skew-symmetric map $f$ corresponds to an element
$\varphi\in \wedge^2F$, such that for all $a^\ast\in F^\ast$ we have
$f(a^\ast)=-a^\ast(\varphi)$.  
In terms of a basis $e_1,\ldots,e_{2m}$ of $F$, if $(f_{ij})$ is the
matrix of $f$ with $f_{ij}=-f_{ji}$, then
$$\varphi=\sum_{i<j} f_{ij} e_i\wedge e_j.$$

Now fix an orientation $e^\ast\in\wedge^{2m}F^\ast$ (that is
a generator of this module). This yields a correspondence between
skew-symmetric maps $g:F\rightarrow F^*$ and elements $\psi\in \bigwedge^{2m-2}
F$, via $g(a)=\psi(a(e^*))$. In coordinates, write 
$$\psi=\sum_{i<j} g_{ij} e_{\{i,j\}^\ast},$$
where $\{i,j\}^\ast$ denotes the complement of the set $\{i,j\}$ in
$\{1,\ldots,2m\}$, and for a subset $H=\{i_1,\ldots,i_n\}\subseteq
\{1,\ldots,2m\}$ with $i_1<\cdots<i_n$, $e_H=e_{i_1}\wedge\cdots\wedge
e_{i_n}$. Notice that the matrix corresponding to $g$ 
is $((-1)^{i+j}g_{ij})$. We then have,
for $a\in F$,
$$\eqalign{(f\circ g)(a)&=f(\psi(a(e^*)))\cr
&=-\psi(a(e^*))(\varphi).\cr}$$

Recall that divided powers are related to Pfaffians by the formula
$$\varphi^{(p)}=\sum_{|H|=2p}\Pf(M_H)e_H,$$
where $M$ denotes
the skew-symmetric matrix of $f$ with respect to  $\{e_i\}$
and the dual basis $\{e^\ast_i\}$, while $M_H$ denotes the principal
submatrix of $M$ determined by rows and columns indexed by $H$.
Thus, if $M$ has rank $2m-2$, then $\varphi^{(m)}=0$. Now 
Lemma 2.4 of [BE] tells us that for $a\in F$, 
$$\varphi^{(m)}(a(e^*))=\varphi^{(m-1)}(a(e^*))(\varphi),$$
so if $g:F\rightarrow F^*$ is the morphism corresponding to $\varphi^{(m-1)}$
and $\varphi^{(m)}=0$, then $f\circ g=0$.

In coordinates, let $M$ be the skew symmetric $n\times n$-matrix
corresponding to $f$ (and the above choice of bases),
denote by $\Pf(M)$ its Pfaffian and by $\Pf^{ij}(M)$, for $i<j$, the 
Pfaffian of the skew-symmetric matrix obtained from
$M$ by ommiting the $i^{\rm th}$ and $j^{\rm th}$ rows and columns, 
and set $\Pf^{ij}(M)=-\Pf^{ji}(M)$, if $i>j$.

Thus if the matrix $M^\ast=(m^\ast_{ij})$  is 
the  skew-symmetric matrix with entries
$$m^\ast_{ij}=
\cases{{(-1)}^{i+j}\Pf^{ij}(M)& if\ $i\ne j$,\cr
0& if\ $i=j$,\cr}$$
then  $M^\ast$ is the matrix of $g$ 
associated with the bases $\{e_i\}$ of $F$
and $\{e_i^\ast\}$ of $F^\ast$. The above compositions read as:
$M\cdot M^\ast=0$. Thus also $M^*\cdot M=(M\cdot M^*)^t=0$, as
required.\Box

We show next that $\overline{\im(\Theta)}$ is a smooth Fano 3-fold of 
genus 8 and index one. As a corollary of \ref{rat.map11} and 
\ref{lin.eq} below we obtain then \ref{main}.

\theorem{lin.eq} The Zariski closure of 
$\im(\Theta)\subset\Gr(2,6)=\Gr(2,V_+)$ in the
Pl\"ucker embedding has equations given in Pl\"ucker coordinates by
$$
p_{23}=-p_{15}, \qquad
p_{26}=p_{13}, \qquad
p_{14}=-p_{35},\qquad
p_{16}=p_{45},\qquad
p_{46}=-p_{12}.$$
Furthermore, this linear section of $\Gr(2,6)$ is three-dimensional,
smooth, and hence a Fano 3-fold of type $V_{14}$. Furthermore it
is birational to the Klein cubic hypersurface
$$\K=V(\sum_{i\in\boldz_5} x_i^2x_{i+1}=0)\subset\Pfour.$$ 

\proof: The first thing to check is that $\im(\Theta)\subset\P^9$ satisfies the
five linear relations given above. This can easily be checked
by hand using Lemma 2.5
simply by computing the corresponding Pfaffians of the matrix $S$
and showing they satisfy the given relations. Observe also that the 
$\SL_2(\boldz_{11})$ representation on $\wedge^2(V_+)$ induced by 
$\rho_+$ decomposes as the sum of a 5 and a 10-dimensional 
irreducible representations. The equations in the 
statement of \ref{lin.eq} define this 10-dimensional 
representation as a subspace of
$\wedge^2(V_+)$. 
 
To conclude that the closure of the image of $\Theta$ is actually
given by these equations, we will show that the subscheme $X$ of $\Gr(2,6)$
defined by these equations is three-dimensional and non-singular, and thus
a Fano 3-fold $V_{14}$ of genus 8 and index 1. To this end we will use a classical construction due to 
G. Fano [Fa], and recast in modern language
by Iskovskih [Is1], [Is2] (see also [Pu]), which
shows that any $V_{14}$ is birationally equivalent to a (smooth) cubic threefold in $\P^4$. 

Let $H_1,\ldots,H_5$ be five (linearly independent) hyperplanes in $\P^{14}=\P(H^0(\O_{\Gr(2,6)}(1)))$,
and let $$X:=\Gr(2,6)\cap H_1\cap\cdots\cap H_5.$$ By \ref{chords} and the discussion
preceding it, we may identify $\dual{(\P^{14})}$ with the space of $6\times 6$
skew-symmetric matrices, and so $\dual{\Gr(2,6)}$ can be naturally
identified with the locus $R_2$ of skew-symmetric matrices of rank $\le 4$.
As seen above $\dual{\Gr(2,6)}$ is then a cubic hypersurface in $\dual{(\P^{14})}$
defined by the $6\times 6$ Pfaffian of the generic skew-symmetric matrix.
Its singular locus is the locus $R_1$ of $6\times 6$
skew-symmetric matrices of rank $\le 2$, which is isomorphic 
(as an embedded variety) to $\Gr(2,6)$. By \ref{chords}, the Pfaffian
cubic is also the secant variety to  $R_1\cong\Gr(2,6)$.

Let $\P^4:=\langle H_1,\ldots,H_5\rangle \subseteq \dual{(\P^{14})}$ denote the span of 
the above five hyperplanes as points in $\dual{(\P^{14})}$, and let
$$B:=\dual{\Gr(2,6)} \cap \langle H_1,\ldots, H_5 \rangle.$$
Define now a (possibly rational) map  
$$\Psi :B\rTo \Gr(2,6)$$ 
$$B\ni b\rMapsto \ker(b)\in \Gr(2,6),$$ 
where we think of each element $b$ of $B$ as a
$6\times 6$ skew-symmetric matrix of rank $\le 4$.
Therefore $\Psi$ is defined on all of $B$ iff $B$ is disjoint
from the singular locus of $\dual{\Gr(2,6)}$.

\lemma{smooth} If the cubic hypersurface $B$ is smooth and $\im(\Psi)\cap X=\emptyset$, 
then $X$ is a non-singular threefold (Fano of genus 8, index 1).

\proof: First, let $x\in X$ be a point where the Zariski tangent space of
$X$ at $x$, $T_{X,x}$, has $\dim T_{X,x}>3$. If $l_1,\ldots,l_5$ are the
Zariski tangent spaces to $H_1,\ldots, H_5$ at $x$, then 
$$T_{X,x}=T_{\Gr(2,6),x}\cap l_1\cap \ldots\cap l_5\subseteq T_{
\P^{14},x}.$$
Suppose now that $H_1,\ldots, H_5$ are given by the linear 
equations $h_1=0,\ldots, h_5=0$, respectively. 
The only way $T_{X,x}$ could fail to be
three dimensional is if there exist a hyperplane section $H$ whose equation
is $\sum_{i=1}^5 a_i h_i=0$ for some $a_i$, 
with $T_{H,x}\supseteq T_{\Gr(2,6),x}$. Thus 
$H$ must be tangent to $\Gr(2,6)$, 
and so $H\in B$. If the cubic $B$ is smooth, 
then $B$ is disjoint from the singular 
locus of $\dual{\Gr(2,6)}$, and the map $\Psi$ above is defined 
everywhere on $B$. Now $\Psi(H)$ is the point of the 
Grassmannian $\Gr(2,6)$ which $H$ is tangent to; 
thus, in particular, if $H$ is tangent to $\Gr(2,6)$ at a point of $X$,
we must have $\Psi(H)\in X$ and $\Psi(B)\cap X\not=\emptyset$. 
In conclusion, if $\Psi(B)\cap X=\emptyset$, we must have  
$\dim T_{X,x}=3$ and so $X$ is a
non-singular Fano threefold. \Box

{\it Proof of \ref{lin.eq} continued:}  We use \ref{smooth} to check 
that the subscheme $X$ of $\Gr(2,6)$
defined by the equations in the statement of 
\ref{lin.eq} is three-dimensional and non-singular.
We can now compute $B$ directly in our case: Let $x_{ij}$, $1\le i<j \le 6$,
be coordinates on $\dual{(\P^{14})}$ dual to the Pl\"ucker coordinates 
$p_{ij}$ on $\P^{14}$. Then, in our particular case, 
the $\Pfour$ spanned by $H_1,\ldots,H_5$ 
is cut out by the equations
$$\eqalign{ x_{12}-x_{46}&=0,\cr
x_{13}+x_{26}&=0,\cr
x_{14}-x_{35}&=0,\cr
x_{15}-x_{23}&=0,\cr
x_{16}+x_{45}&=0,\cr
x_{24}=x_{25}&=x_{34}=x_{36}=x_{56}=0.\cr}$$
Making now the substitutions $x_0=x_{12}, x_2=x_{13}, x_1=x_{14}, 
x_4=x_{15},$ and $x_3=x_{16}$, we see that the equation of $B\subseteq \Pfour$,
with coordinates $x_0,\ldots, x_4$, is given by the $6\times 6$ Pfaffian
of the skew-symmetric matrix
$$M=\pmatrix{ 0& x_0& x_2 & x_1 & x_4 & x_3\cr
-x_0 & 0 & x_4 & 0 & 0 & -x_2\cr
-x_2 & -x_4 & 0 & 0 & x_1 & 0\cr
-x_1 & 0 & 0 & 0 & -x_3 & x_0\cr
-x_4 & 0 & -x_1 & x_3 & 0 & 0\cr
-x_3 & x_2 & 0 & -x_0 & 0 & 0\cr
}$$
which is
$$B=\{x_0^2 x_1+ x_1^2 x_2 + x_2^2 x_3 +x_3^2 x_4 + x_4^2 x_0=0\}.$$
Thus $B$ is Klein's cubic $\K=\{\sum_{i=0}^4 x_i^2 x_{i+1}=0\}$, the only
$PSL_2(\boldz_{11})$-invariant cubic in $\Pfour$. 
This cubic is known to be smooth. (See also [Ad4], Lemma 47.2 
for the Pfaffian description of the Klein cubic.) 

To show that $X$ is non-singular, we need now to check the second
hypothesis of \ref{smooth}. By \ref{pfaff}, 
the map $\Psi:B\rightarrow \Gr(2,6)$ is given in Pl\"ucker
coordinates by the matrix
$$(M)^*=
\pmatrix{0&x_0x_1&x_2x_3&x_1x_2&x_0x_4&x_3x_4\cr
-x_0x_1&0&x_3^2+x_0x_4&x_1x_3&-x_0x_2&-x_1^2-x_2x_3\cr
-x_2x_3&-x_3^2-x_0x_4&0&-x_2x_4&x_0^2+x_1x_2&x_0x_3\cr
-x_1x_2&-x_1x_3&x_2x_4&0&-x_2^2-x_3x_4&x_0x_1+x_4^2\cr
-x_0x_4&x_0x_2&-x_0^2-x_1x_2&x_2^2+x_3x_4&0&-x_1x_4\cr
-x_3x_4&x_1^2+x_2x_3&-x_0x_3&-x_0x_1-x_4^2&x_1x_4&0\cr}.$$
In order for a point $P=(x_0:\ldots:x_4)\in B$ to satisfy $\Psi(B)\in X$,
the Pl\"ucker coordinates of $\Psi(P)$ must satisfy the five linear
equations defining $X$, which yields that $P$ must satisfy the equations
$$\eqalign{x_i^2+2x_{i+1}x_{i+2}&=0,\quad\quad 0\le i \le 4\cr
\sum_{i=0}^4 x_i^2x_{i+1}&=0\cr}.$$
The first set of equations is precisely the Jacobian of the Klein cubic
$B$, and since $B$ is smooth, there are no points $P\in B$ satisfying these
equations. Hence, by \ref{smooth}, $X$ is non-singular. 

We are left now to construct a birational map between $B(=\K)\subset\P^4$ and 
$X\subset\P^9$. As mentioned above, a classical construction due to Fano and 
Iskovskih [Fa], [Is2]
(see also [Pu] for details) provides such a (rather indirect)  birational transformation.   
For the reader's convenience we sketch it in the sequel.

The lines $L_p$ in $\P^5=\P(\dual{V_+})$ 
represented by points $p\in X$ sweep out
an irreducible quartic hypersurface $\Gamma\subset\P^5$ (called ``da Palatini''
by Fano). Through the generic point 
$x\in\Gamma$, passes exactly one $L_p$, with $p\in X$.

On the other hand, since $B$ is smooth, each point $q\in B=\dual{\Gr(2,6)}\cap\P^4$ 
corresponds to a hyperplane $H_q$ tangent to the Grassmannian in exactly one
point $n_q$, called the ``centre'' of $H_q$. 
The lines $N_q$ in $\P^5$ represented
by centres $n_q$ of points $q\in B$ sweep out 
an irreducible variety $\Sigma\subset\P^5$.

It is easy to see that 
$\Sigma\subset \Gamma$. (See for instance [Pu], pp. 83--84,
where the given argument holds whenever $B$ is a smooth 3-fold.) 
Now through the generic
point of $\Sigma$ passes exactly one line $N_q$, with $q\in B$. Otherwise, if
$q, q'\in B$ and  $N_q\cap N_{q'}\ne\emptyset$, then the whole pencil
spanned by $q$ and $q'$ lies in $B$, which means that we found a line in $B$.
But this contradicts the fact that the Fano variety of the  
Klein cubic is 2-dimensional.
Therefore $\Sigma\subset\P^5$ is  an irreducible hypersurface, and so
we must have $\Sigma=\Gamma$. (The quartic equation defining $\Gamma=\Sigma$
is the unique quartic invariant for the action of 
$\SL_2(\boldz_{11})$ on $\P^+$ via $\rho_+$; see [Ad4], Corollary 50.2
for the explicit equation.) 

Choose now a generic hyperplane $\Pi\subset\P^5$, 
and let $\bar\Gamma:=\Gamma\cap\Pi$.   
We may define birational maps
$$\eta: X \rDashto \bar\Gamma,\qquad \eta(p):=\Pi\cap L_p,$$
$$\gamma: B \rDashto \bar\Gamma,\qquad \gamma(q):=\Pi\cap N_q.$$
The composition $\chi:=\gamma^{-1}\circ\eta$ defines now a birational
isomorphism between $X$ and $B$, as required. (See [Pu]
and [Is2] for a detailed analysis of this mapping.) \Box

\remark{final} 1) The birational isomorphism $\chi$ provided in the proof of 
\ref{lin.eq} and \ref{main} depends on the choice of a hyperplane $\Pi\subset\P^5$,
and thus is not compatible with the action of $\PSL_2(\boldz_{11})$. 
The indeterminacy locus of the isomorphism $\chi$ (as well as of its
inverse) turns out to be the union of an elliptic quintic curve $E$
and 25 mutually disjoint secant lines to it (which are flopped by $\chi$).
In terms of linear systems, $\chi$ is defined by $|5H-3E|$, where
$H$ is the hyperplane class on $X$. Similarly, $\chi^{-1}$
is induced by $|7H'-4E'|$, where $H'$ is the hyperplane class on $B$
and $E'$ is the base locus of $\gamma$.\hfill\break
2) Takeuki [Ta] and Tregub [Tr] have constructed a different
birational isomorphism of a smooth Fano 3-fold $V_{14}$
of genus 8, index one  onto a smooth cubic hypersurface $B\subset\P^4$, which
can be briefly described as follows. Let $C$ be a (general) rational normal 
curve  $C$ on $B\subset\P^4$. There are exactly 16 chords $l_i, 1\le i\le 16$ to $C$ on $B$. 
Let $\widetilde B$ be the blowing-up of $B$ along $C$ and the $l_i$'s. 
Then the linear system $L=|8H-5C-\sum_{i=1}^{16}2L_i|$ provides a birational 
morphism from $\widetilde B$ onto the intersection $X$ of $\Gr(2,6)$ with a codimension 5 linear
subspace. Under this morphism, the unique divisor $D\in|3H-2C-\sum_{i=1}^{16}L_i|$ 
is contracted to a point $p$. The inverse
birational morphism is then induced by the linear system 
$\vert 2H'-3p|$, where $H'$ is the hyperplane class on $X$.\hfill\break 
3) Notice also that every (abstract) smooth Fano 3-fold  of genus 8, 
index one is isomorphic to a codimension 5 linear
section of  $\Gr(2,6)$, cf. [Gu].

\question{jacobian} It seems plausible that the lines on the
codimension 5 linear section $Y$ of  $\Gr(2,6)$, which is the Zariski closure
of $\im(\Theta)$ in \ref{lin.eq}, are parametrized by the
modular curve $X(11)$, and that the intermediate Jacobian
of $Y$ is isomorphic to the generalized Prym variety
corresponding to the  (symmetric) Hecke
correspondence $T_3$ on $X(11)$. See, for instance,  [Ad2] and [Ed]
for a geometric description of this Hecke correspondence.

\remark{kummer11}  Let $(A,\L)$ be a general
$(1,11)$-polarized abelian surface, where $\L$ is assumed to be
symmetric. One can  show that the linear system 
${|2\L-2\sum_{i=1}^{16} e_i|}^+$,
of even divisors of  $2\L$ having
multiplicity two in the half periods,
descends to a (complete) very ample linear
system on the (desingularized) Kummer surface $X$ associated to $A$,
and embeds it as a codimension 8 linear section of the 
spinor variety ${\cal S}\subset\P^{15}$,
which parametrizes isotropic $\Pfour$'s in an 8-dimensional
smooth quadric in $\P^9$. It would be interesting, in the light of
[Muk], to determine exactly which codimension 8 linear sections correspond
to such Kummer surfaces.

\section {moduli9} {Moduli of $(1,9)$-polarized abelian surfaces}

As mentioned in the introduction, 
an argument similar to the one used in \ref{moduli11} 
allows us to prove the rationality of $\A_9^{lev}$.

For the remaining of the paper, let $(A,\L)$ be a general abelian surface 
with a polarization of type $(1,9)$ and with canonical level structure. 
Most of the facts concerning theta groups from \ref{representations}
can be adapted to this case, but we will make little use of them in the sequel.
We will also assume that $\L$ is chosen 
to be a symmetric line bundle.

As seen in \ref{representations} 
$|\L|$ embeds $A\subset\P^{8}=\P(\dual{H^0(\L)})=\P(\dual{V})$
invariantly under the action of the
Heisenberg group $\HHH_{9}$ via the Schr\"odinger representation,
and the involution $\iota$. In particular, $H^0(\I_A(n))$ is a 
representation of weight $n$ of the Heisenberg group, whose 
irreducible components will have dimension $9/\gcd(9,n)$.
Via $|\L|$, $A$ is embedded as a projectively normal surface
of degree 18  which is  contained in $9$ quadrics 
(cf. [Laz], or [GP1], Theorem 6.5). However, in contrast with \ref{moduli11},
we are in a boundary case, in that these quadrics do not generate 
the homogeneous ideal of $A$.  Moreover, 
in general the quadrics containing the degenerations
used in the proof of \ref{rat.map11} and [GP1], Theorem 6.5, b),  
cut out only a threefold.

As in \ref{moduli11} we will investigate the locus of odd 2-torsion points,
which in this simpler case turns out to be the whole of 
$\P^-=\P(\dual{V_-})\cong\P^3$. 

The space of quadrics $H^0(\O_{\P^8}(2))$ decomposes into five
9-dimensional representations of the Heisenberg group, each one 
isomorphic to the Schr\"odinger representation. As above, 
one such decomposition is given
by the spans of the rows of the matrix defined in [GP], \S 6:
$$R_4=\pmatrix{
x_0^2&x_1^2&x_2^2&x_3^2&x_4^2&x_5^2&x_6^2&x_7^2&x_8^2\cr
x_1x_8&x_2x_0&x_3x_1&x_4x_2&x_5x_3&x_6x_4&x_7x_5&x_8x_6&x_0x_7\cr
x_2x_7&x_3x_8&x_4x_0&x_5x_1&x_6x_2&x_7x_3&x_8x_4&x_0x_5&x_1x_6\cr
x_3x_6&x_4x_7&x_5x_8&x_6x_0&x_7x_1&x_8x_2&x_0x_3&x_1x_4&x_2x_5\cr
x_4x_5&x_5x_6&x_6x_7&x_7x_8&x_8x_0&x_0x_1&x_1x_2&x_2x_3&x_3x_4\cr
}.$$
Thus every 9-dimensional $\HHH_9$-subrepresentation of quadrics 
is spanned by $v\cdot R_4$ 
for some $v\in{\bf C}^5=V_+$, and thus these representations 
are parametrized by $\P^+:=\P(\dual{V_+})$. 
If we restrict $R_4$ to $\P^-=\P(\dual{V_-})$, the $(-1)$-eigenspace 
of the involution $\iota$, and consider as 
before the first $5\times 5$ block, we obtain 
the matrix
$$S=\pmatrix{
0&x_1^2&x_2^2&x_3^2&x_4^2\cr
-x_1^2&0&x_3x_1&x_4x_2&-x_4x_3\cr
-x_2^2&-x_3x_1&0&-x_4x_1&-x_3x_2\cr
-x_3^2&-x_4x_2&x_4x_1&0&-x_2x_1\cr
-x_4^2&x_3x_4&x_2x_3&x_1x_2&0\cr
},$$
representing the intertwining 
operator $\Phi : \wedge^2(V_+)\rightarrow S^2(V_-).$ 

As in \ref{moduli11}, or [GP1], \S 6, it follows that for a point 
$P\in\P^-$, $v\cdot R_4(P)=0$ if and only if $v\cdot S(P)=0$. 

\lemma{modular9} 
\item{1)} $\rank(S(P))\ge 2$, for all $P\in\P^-$.
\item{2)} The locus  $D_1\subseteq\P^-$ where 
matrix $S$ has rank $2$ is the disjoint union
of a smooth curve $C\subset\P^-$ of degree 9, 
which is the complete intersection
$$\{x_1^2x_2-x_2^2x_4-x_1x_4^2=x_1x_2^2-x_3^3+x_1^2x_4-x_2x_4^2=0\}\subset\P^-,$$ 
and the four points 
$$P_1=(0:0:0:1:0:0:-1:0:0),\qquad P_2=(0:-1:1:0:-1:1:0:-1:1),$$
$$P_3=(0:-1:\xi^3:0:-\xi^6:\xi^6:0:-\xi^3:1),\ P_4=(0:-1:\xi^6:0:-\xi^3:\xi^3:0:-\xi^6:1),$$
where $\xi=e^{{2\pi i}\over 9}$ is a primitive root of order nine of the unity.
The curve $C$ is isomorphic to the modular curve $X(9)$.

\proof: Direct computation and  arguments similar to those
used in the proof of \ref{modular}.\Box

As in \ref{moduli11}, we may interpret $D_1$ as the locus
of points in $\P^-\subset\P^{8}=\P(\dual{V})$ 
which are contained in a net of
$\HHH_{9}$-subrepresentations of quadrics. On the other hand,
$S$ is a $5\times 5$ skew-symmetric matrix, so $S$ drops rank on all of
$\P^-$.  Therefore we can define again a map 
$$\eqalign{\Theta : \P^-\setminus D_1&\rTo \P^+\cr
\P^-\setminus D_1\ni P &\rMapsto \P(\ker(S(P)))=
\P(\{v\in V_+ \mid v\cdot S(P)=0\}),\cr}$$ 
which sends a point $P$ to the unique $\HHH_{9}$-subrepresentation
of $H^0(\O_{\P^8}(2))$ of quadrics containing it. 

By an argument similar to the one used in the proof of \ref{pfaff}, the
morphism $\Theta$ is easily seen to map a point  $P\in\P^-$ 
to the point of $\P^+$ whose coordinates are 
given by the $4\times 4$-Pfaffians of 
$S(P)$, taken with suitable signs. In coordinates, if  
$P=(x_1:\ldots:x_4)\in\P^-$ this yields
$\Theta(x_1:\ldots:x_4)=(v_0:\ldots:v_4)$, where
$$
\eqalign{
v_0&=-x_1^2x_2x_3+x_2^2x_3x_4+x_1x_3x_4^2\cr
v_1&=x_1x_2^3-x_2x_3^3+x_1x_4^3\cr
v_2&=-x_1^3x_2+x_3^3x_4+x_2x_4^3\cr
v_3&=x_1^2x_2x_3-x_2^2x_3x_4-x_1x_3x_4^2\cr
v_4&=x_1x_3^3-x_1^3x_4-x_2^3x_4\cr}.
$$

Since $v_0=-v_3$, we deduce that the image of $\Theta$ is contained in the
linear subspace $\Pi$ of $\P^+$ defined by $v_0=-v_3$. As in \ref{moduli11},
by [GP1], Lemma 6.4, $\Theta$ induces a rational map
$$\Theta_9:\A_9^{lev}\rDashto\Pi,$$
which essentially is defined by taking an abelian surface $A\subseteq\P^8$ to
$\Theta((A\cap\P^-)\setminus D_1)$, that is to the point corresponding 
to the unique $\HHH_{9}$-subrepresentation of $H^0(\O_{\P^8}(2))$ of 
quadrics containing the abelian surface. 

\remark{quadrics9} It is easy to see that a $(1,9)$-polarized abelian surface 
$A\subseteq\P^8$ is not cut out by quadrics. 
Indeed, if $v=(v_0:\ldots:v_4)=\Theta_9(A)\in \im(\Theta)$, then $v_0=-v_3$ 
and each quadric entry of $v\cdot R_4$ vanishes at the (fixed) point
$$P=(1:0:0:1:0:0:1:0:0)\in\P^8.$$
However, since $\sigma^3(P)=P$ and $\tau^3(P)=P$, where $\sigma$ and
$\tau$ acting by translation by 9-torsion points on $A$ are 
the usual generators of 
$\HHH_{9}$ in the Schr\"odinger representation, 
$P$ cannot be contained in $A$. Thus $A$ is not cut
out by quadrics since the only quadrics containing $A$ are linear
combinations of the entries of $v\cdot R_4$. In fact one may show that
for the general abelian surface $A$, the quadrics defined by
$v\cdot R_4$ cut out the union of $A$ and the set of nine points which form
the $\HHH_9$ orbit of $P$. A degeneration argument, in the spirit
of [GP1], \S 6, shows that the homogeneous ideal of $A$ is in fact
generated by the 9 quadrics and 6 extra cubics (use \ref{torus} below).

We can now prove the rationality of $\A_9^{lev}$:

\theorem{rat.map9} 
$\Theta_9:\A_9^{lev}\rDashto\Pi\cong\P^3$ is a birational
map.

\proof: We follow much the same strategy as the proof of [GP1], Theorem 6.5,
however since, by \ref{quadrics9}, quadrics do not cut out an abelian surface, 
we will need to involve cubic equations, and the process is a bit
more difficult computationally. 

We will also make use of the ubiquitous (Moore) $9\times 9$-matrices 
$$M'_4(x,y)={(x_{5(i+j)}y_{5(i-j)})}_{i,j\in\boldz_9},$$
where we think of $x={(x_i)}_{i\in\boldz_9}$ as a point in 
the ambient $\P^8$ and $y={(y_i)}_{i\in\boldz_9}$ as a parameter point.
We refer the reader to [GP1], \S 2 and \S 6 for a detailed discussion
of their properties. Note also that the matrix $R_4$ above, 
up to transpose and 
permutations of rows and columns, is a submatrix of $M'_4(x,x)$.

Let $Z:=\Theta^{-1}(\im(\Theta_9))\subseteq \P^-\setminus D_1,$ and let
$\bar Z$ denote  the closure of $Z$ in $\P^-\setminus D_1$. Let $\A\subseteq
\P^8_{\bar Z}$ be the family defined by the condition that the ideal of 
a fibre $\A_z$, $z\in\bar Z$, is generated by the $9$ quadrics, which are
entries of $\Theta(z)\cdot R_4$ (i.e., the $\HHH_9$-subrepresentation of
$H^0(\O_{\P^8}(2))$ vanishing at $z$), along with the cubics which are the 
$6\times 6$-Pfaffians of the skew-symmetric $9\times 9$-matrix
$M_4'(x,z)$. By [GP1], Corollary 2.8 and Lemma 6.4, $\A_z$ contains all abelian
surfaces whose odd 2-torsion points map to $z$.

We need to show that there exists an open set $U\subset \bar Z$, such that
the restricted family $\A_U\rTo U$ is flat, and every smooth
fiber is an $\HHH_9$-invariant (and thus $(1,9)$-polarized) abelian surface.

The degeneration argument in [GP1], Theorem 3.1 and Theorem 6.2 shows that if 
$E\subset\P^8$ is a Heisenberg invariant elliptic normal
curve of degree 9, then $\Sec(E)\cap (\P^-\setminus D_1)\subseteq \bar{Z}$.
The same is also true if we take $E$ to be the 
``standard 9-gon'' $X(\Gamma_9)$,
and $\Sec(E)$ to be its ``secant variety'', 
that is, with notation as in [GP1], \S 4: 
$$X(\Gamma_9)=\cup_{i\in\boldz_{9}} l_{i,i+1}\subseteq\P^8,$$
where  $l_{i,i+1}=\langle e_i,e_{i+1}\rangle$
is the line joining the vertices $e_i$ and $e_{i+1}$ of
the standard simplex in $\P^8$. In particular, 
for $E=X(\Gamma_9)$, the set
$\Sec(E)\cap(\P^-\setminus D_1)$ and thus also $\bar Z$ contain 
the point $$z_0=(0:0:-1:-1:0:0:1:1:0).$$

Let $I_0$ be the homogeneous ideal of the fibre $\A_{z_0}$. 
To conclude the result, it will be enough to show that $\A_{z_0}$ 
is contained in a surface of degree $18$. 

Now $\Theta(z_0)=(0:1:0:0:0)\in\Pi\subset\P^+$, 
so $I_0$ contains the quadrics
$\Theta(z)\cdot R_4$, namely:
$$\{x_ix_{i+2}, \quad i\in\boldz_9\}.$$
On the other hand, the matrix $M'_9(x,z_0)$ is
$$\pmatrix{0&0&0&{-{{x}}_{{6}}}&{{x}}_{{2}}&{-{{x}}_{{7}}}&{{x}}_{{3}}&0&0\cr 
0&0&0&0&{-{{x}}_{{7}}}&{{x}}_{{3}}&{-{{x}}_{{8}}}&{{x}}_{{4}}&0\cr 
0&0&0&0&0&{-{{x}}_{{8}}}&{{x}}_{{4}}&{-{{x}}_{0}}&{{x}}_{{5}}\cr 
{{x}}_{{6}}&0&0&0&0&0&{-{{x}}_{0}}&{{x}}_{{5}}&{-{{x}}_{1}}\cr 
{-{{x}}_{{2}}}&{{x}}_{{7}}&0&0&0&0&0&{-{{x}}_{1}}&{{x}}_{{6}}\cr 
{{x}}_{{7}}&{-{{x}}_{{3}}}&{{x}}_{{8}}&0&0&0&0&0&{-{{x}}_{{2}}}\cr 
{-{{x}}_{{3}}}&{{x}}_{{8}}&{-{{x}}_{{4}}}&{{x}}_{0}&0&0&0&0&0\cr 
0&{-{{x}}_{{4}}}&{{x}}_{0}&{-{{x}}_{{5}}}&{{x}}_{1}&0&0&0&0\cr 
0&0&{-{{x}}_{{5}}}&{{x}}_{1}&{-{{x}}_{{6}}}&{{x}}_{{2}}&0&0&0\cr }$$
We consider first two of its $6\times 6$-Pfaffians: 
The skew-symmetric $6\times 6$-minor 
coming from taking rows (and columns) 1,2,3,5,6 and 7 has Pfaffian
$$-x_2x_3x_4+x_4x_7^2-x_3x_7x_8+x_2x_8^2\in I_0,$$
and similarly taking rows (and columns) 1,2,3,4,6, and 8, we get
another cubic Pfaffian
$$-x_0x_3x_6+x_4x_6x_8\in I_0.$$
Taking into account the quadrics in $I_0$, we observe that $I_0$ also contains the polynomials
$x_4x_7^2-x_3x_7x_8+x_2x_8^2$ and $x_0x_3x_6$. 
Since the matrix $M_9'$ is Heisenberg
invariant (in the $x$-coordinate) up to permutations of rows and columns,
it follows that $I_0$ is $\HHH_9$-invariant, and hence 
contains
$$\eqalign{x_ix_{i+2},&\quad i\in\boldz_9\cr
x_0x_3x_6,\ \ x_1x_4x_7,\ \ x_2x_5x_8,&\cr
x_{i+4}x_{i+7}^2-x_{i+3}x_{i+7}x_{i+8}+x_{i+2}x_{i+8}^2,&
\quad i\in\boldz_9.\cr}$$ 
The claim of \ref{rat.map9} follows now from the following
combinatorial lemma, which determines the Hilbert polynomial of $I_0$:

\lemma{torus}
\item{1)} The ideal $J_1$ generated by the quadric and cubic monomials 
$$\{x_ix_{i+2},\quad x_ix_{i+3}x_{i+6},\quad x_{i+3}x_{i+7}x_{i+8}
\mid  i\in\boldz_9\}$$
is the Stanley-Reisner face ideal $I_{X(\Delta_9)}$ corresponding to the 
triangulation $\Delta_9$ of the torus $T_1$ in [GP1], Proposition 4.4. 
In particular, $J_1$ has the same
Hilbert polynomial as a $(1,9)$-polarized abelian surface.
\item{2)}  The ideal $J_2$ generated by the 12 quadric and cubic monomials 
$$\{x_ix_{i+2},\quad x_ix_{i+3}x_{i+6}\mid i\in\boldz_9\},$$
cuts out the threefold 
$$\Sigma=\bigcup_{i=0}^8 L_i,\quad\quad L_i=\sigma^i(L_0),$$
where $L_0$ is the $\Pthree$ determined by $\{x_0=x_1=x_4=x_5=x_8=0\}.$
$J_2$ is the face ideal of the ``solid'' torus whose triangulation
$\Delta_9$ is described in $1)$. This is the complex whose two-simplices
are those of $\Delta_9$ but which has in addition three-simplices 
with vertices $(x_i,x_{i+1},x_{i+3},x_{i+4})$.
\item{3)} The ideals
$$J_{(\lambda:\mu)}=J_2+\langle  
\lambda x_{i+4}x_{i+7}^2-\mu x_{i+3}x_{i+7}x_{i+8}+
\lambda x_{i+2}x_{i+8}^2,\ i\in\boldz_9
\rangle,$$ 
for $(\lambda:\mu)\in\P^1$, define a flat family of surfaces 
$X_{(\lambda:\mu)}\subset\P^8$
with  the same Hilbert polynomial as a $(1,9)$-polarized abelian surface.
In particular, $I_0=J_{(1:1)}$ defines a surface of degree 18 as desired.

\proof: The proof is easy and left to the reader. Observe
that $X_{(\lambda:\mu)}$ is defined by $J_2$ and
9 trinomials, from which it can be shown that
set theoretically $X_{(\lambda:\mu)}$ is the union
of 9 distinct (smooth) quadric surfaces
$$X_{(\lambda:\mu)}=\bigcup_{i=0}^8 Q_i,\quad\quad Q_i=\sigma^i(Q_0),$$
where $Q_0$ is defined by
$$Q_0=L_0\cap \{\lambda x_3x_6-\mu x_2x_7=0\}.$$
On the other hand, $J_{(0:1)}=J_1$ is the 
Stanley-Reisner ideal of the triangulation
$\Delta_9$ of the torus, and thus has the required Hilbert function.
(See also [GP1], Proposition 4.4, for details.)\Box

\remark{kummer9} The linear projection $\pi_-: \P(\dual{V})=
\P^8\rDashto \P^-$ commutes with the involution 
$\iota$ and thus maps a general 
$\HHH_9\rtimes\langle\iota\rangle$ abelian
surface $A\subset\P^8$ to a 6-nodal Kummer quartic surface 
$K\subset\P^3=\P^-$ (whose nodes are the odd $2$-torsion points 
of $A$). The linear system of quadrics through the set $S$ of nodes
of $K$ maps the Kummer surface to a smooth quartic $K'\subset\P^3$.
Such a smooth quartic surface has 16 skew conics, and in fact
any smooth quartic surface in $\P^3$ containing 16 skew conics
is a Kummer surface of an abelian surface $(A,\L)$ with a polarization
of type $(1,9)$, via the linear system ${|2\L-2\sum_{i=1}^{16} e_i|}^+$
of even divisors of the totally symmetric line bundle $2\L$, having
multiplicity two in the half periods (see [BB], Claims 2--4, and 
[Bau], Theorem 2.1 for a detailed discussion).

\remark{pol33} There is a second family of (minimal) abelian surfaces
of degree 18, and sectional genus 10 embedded in $\P^8$, namely
those embedded via a polarization of type $(3,3)$. These are also
contained in 9 independent quadrics, that, in contrast with
the $(1,9)$ case, cut out scheme theoretically the abelian surface.
The homogeneous ideal of a $(3,3)$ polarized abelian surface
is generated by (quadrics and three independent) cubics (cf. [Se]).
See [Co], [Gra], [vdG], and [Ba] for explicit equations and their relation
to the Burchardt quartic.

\bigskip
\references

\item{[Ad1]}
Adler, A., ``On the automorphism group of a certain cubic threefold'',
{\it Amer. J. Math.},  {\bf 100}, (1978), no. 6, 1275--1280.
\item{[Ad2]}
Adler, A., ``Modular correspondences on $X(11)$'', 
{\it Proc. Edinburgh Math. Soc.}, (2) {\bf 35}, (1992), no. 3, 
427--435.
\item{[Ad3]}
Adler, A., ``Invariants of ${\rm PSL}_2({\bf F}_{11})$ 
acting on ${\bf C}^5$'', {\it Comm. Algebra},
{\bf 20}, (1992), no. 10, 2837--2862.
\item{[Ad4]}
Adler, A., Appendices to
{\it Moduli of Abelian Varieties},
Lect. Notes in Math. {\bf 1644}, Springer Verlag 1996.
\item{[Ad5]} Adler, A., ``Invariants of $\SL_2({\bf F}_q)\cdot
{\rm Aut}({\bf F}_q)$ Acting on ${\bf C}^n$, $q=2n\pm 1$,'' in
{\it The Eightfold way: The beauty of the Klein Quartic}, 175--219,
MSRI Publications {\bf 35}, Cambridge University Press, New York, 1999.
\item{[AR]} Adler, A., Ramanan, S., {\it Moduli of Abelian Varieties},
Lect. Notes in Math. {\bf 1644}, Springer Verlag 1996.
\item{[Ba]}  Barth, W., ``Quadratic equations for level-$3$ abelian 
surfaces'', in {\it Abelian varieties (Egloffstein, 1993)}, 1--18, de Gruyter, 
Berlin, 1995.
\item{[BB]}  Barth, W., Bauer, Th., 
``Smooth quartic surfaces with $352$ conics'', 
{\it Manuscripta Math.}, {\bf 85}, (1994), no. 3-4, 409--417. 
\item{[Bau]} Bauer, Th., ``Quartic surfaces with $16$ skew conics'', 
{\it J. Reine Angew. Math.},  {\bf 464}, (1995), 207--217. 
\item{[BS]} Bayer, D., Stillman, M.,
``Macaulay: A system for computation in
        algebraic geometry and commutative algebra
Source and object code available for Unix and Macintosh
        computers''. Contact the authors, or download from
        {\tt ftp://math.harvard.edu} via anonymous ftp.
\item{[Bea]} Beauville, A.,
``Les singularit\'es du diviseur $\Theta $ de la jacobienne interm\'ediaire de 
l'hypersurface cubique dans $\P^4$'', in {\it Algebraic threefolds (Varenna, 1981)},
190--208,  Lecture Notes in Math., {\bf 947}, Springer Verlag, 1982. 
\item{[BM]} Behr, H., Mennicke, J., ``A presentation of the groups ${\rm PSL}(2,\,p)$'', 
{\it Canad. J. Math.}, {\bf 20}, (1968), 1432--1438.
\item{[Bo]} Borisov, L., ``A fi\-ni\-te\-ness theorem for subgroups
of ${\rm Sp}(4,\boldz)$,'' preprint math.AG/9510002. 
\item{[BE]} Buchsbaum, D., Eisenbud, D.: 
Algebra structures for finite free resolutions, and some structure 
theorems for ideals of codimension $3$. 
{\it Amer. J. Math.} {\bf 99} (1977), no. 3, 447--485.
\item{[CG]} Clemens, C. H., Griffiths, P., 
``The intermediate Jacobian of the cubic threefold'', 
{\it Ann. of Math.}, (2) {\bf 95}, (1972), 281--356. 
\item{[CNPW]} Conway, J. H., Curtis, R. T., Norton, S. P., Parker, R. A.; Wilson, R. A., 
{\it Atlas of finite groups. Maximal subgroups and ordinary characters for
simple groups. With computational assistance from J. G. Thackray}, 
Oxford University Press, Oxford, 1985. 
\item{[Co]} Coble, A., ``Point Sets and Allied Cremona Groups III'',
{\it Trans. Amer. Math. Soc.}, {\bf 18}, (1917), 331--372.
\item{[Dol]} Dolgachev, I., ``Invariant stable bundles over modular
curves $X(p)$'', preprint math.AG/9710012.
\item{[Dor]} Dornhoff, L., {\it Group Representation Theory}, Marcel Dekker, New York 1971/72. 
\item{[Ed]} Edge, W. L., ``Klein's encounter with the simple group of order $660$'', 
{\it Proc. London Math. Soc.} (3), {\bf 24 }, (1972), 647--668.  
\item{[Fa]} Fano, G., Nuove ricerche sulle variet\`a algebriche a tre dimensioni
a curve-sezionicanoniche, {\it Pontificia Acad. Sci. Coment.}, {\bf 11}, (1947), 635-720.
\item{[Gra]}  Grant, D., ``Formal groups in genus two'', {\it J. Reine Angew. Math.}, 
{\bf 411}, (1990), 96--121. 
\item{[GS]} Grayson, D., Stillman, M., 
``Macaulay 2: A computer program designed to support
computations in algebraic geometry and computer algebra.''
Source and object code available from
{\tt http://www.math.uiuc.edu/Macaulay2/}.
\item{[Gri1]} Gritsenko, V., ``Irrationality of the moduli spaces of 
polarized abelian surfaces. With an appendix by the author and K. Hulek''
in {\it Abelian varieties (Egloffstein, 1993)}, 63--84, de Gruyter, Berlin, 1995.
\item{[Gri2]} Gritsenko, V., ``Irrationality of the moduli spaces of 
polarized abelian surfaces'', {\it Internat. Math. Res. Notices} (1994), 
no. 6, 235 ff., approx. 9 pp. (electronic).  
\item{[GH]} Gritsenko, V., Hulek, K., ``Commutator coverings of Siegel threefolds'', 
{\it Duke Math. J.},  {\bf 94}, (1998), no. 3, 509--542. 
\item{[GP1]} Gross,~M., Popescu,~S., 
``Equations of $(1,d)$-polarized abelian surfaces'',
{\it Math. Ann}, {\bf 310}, (1998), (2), 333-377.
\item{[GP2]} Gross,~M., Popescu,~S., ``Calabi-Yau 3-folds and moduli of
abelian surfaces I'', in preparation.
\item{[GP3]} Gross,~M., Popescu,~S., ``Calabi-Yau 3-folds and moduli of
abelian surfaces II'', in preparation.
\item{[Gu]} Gushel, N. P., ``Fano $3$-folds of genus $8$'', 
{\it Algebra i Analiz}, {\bf 4}, (1992), no. 1, 120--134. 
\item{[HM]} Horrocks, G., Mumford, D., `` A Rank 2 vector bundle on  $\Pfour$
with 15,000 symmetries'', {\it Topology}, {\bf 12}, (1973) 63--81.
\item{[Hul]} Hulek, K., {\it Projective Geometry of Elliptic Curves},
{\it Ast\'erisque} {\bf 137}, 1986.
\item{[HKW]} Hulek, K., Kahn, C., Weintraub, S., {\it Moduli Spaces of Abelian
Surfaces: Compactification, Degenerations, and Theta Functions},
Walter de Gruyter 1993.
\item{[HS1]} Hulek, K., Sankaran, G. K., 
``The Kodaira dimension of certain moduli spaces of abelian surfaces'', 
{\it Compositio Math.}, {\bf 90}, (1994), no. 1, 1--35.
\item{[HS2]} Hulek, K., Sankaran, G. K., ``The geometry of Siegel modular varieties'',
preprint math.AG/9810153.
\item{[Is1]} Iskovskih, V. A., ``Fano threefolds. II'', {\it Izv. Akad. Nauk SSSR Ser. Mat.}, 
{\bf 42}, (1978), no. 3, 506--549. 
\item{[Is2]} Iskovskih, V. A., ``Birational automorphisms of three-dimensional algebraic varieties'', in
{\it Current problems in mathematics}, Vol. 12, pp. 159--236, 239, VINITI, Moscow, 1979. 
\item{[Kl]} Klein, F., ``\"Uber die Transformation elfter Ordnung
der elliptische Funktionen, {\it Math. Ann.}, {\bf 15}, (1879),
[Ges. Math. Abh., Band III, art. LXXXVI, 140--168]. 
\item{[KlF]} Klein, F., Fricke, R., {\it Theorie der elliptischen
Modulfunktionen}, Bd. {\bf I}, Teubner, Leipzig 1890.
\item{[Laz]} Lazarsfeld, R., ``Projectivit\'e normale  
des surfaces abeliennes'', (written by O. Debarre), 
preprint Europroj {\bf 14}, Nice.
\item{[LB]} Lange, H., Birkenhake, Ch., {\it Complex abelian varieties},
Springer-Verlag 1992.
\item{[MS]} Manolache, N., Schreyer, F.-O., 
``Moduli of (1,7)-polarized abelian surfaces via syzygies'',
preprint math.AG/9812121.
\item{[Muk]}  Mukai, S., ``Curves, $K3$ surfaces and Fano $3$-folds of genus $\leq 10$'',
in {\it Algebraic geometry and commutative algebra}, Vol. {\bf I}, 
357--377, Kinokuniya, Tokyo, 1988. 
\item{[Mum]} Mumford, D., {\it Abelian varieties}, Oxford 
University Press 1974.
\item{[Mur]} Murre, J. P. 
``Reduction of the proof of the non-rationality of a non-singular 
cubic threefold to a result of Mumford'', 
{\it Compositio Math.}, {\bf 27}, (1973), 63--82. 
\item{[O'G]} O'Grady, K. G., ``On the Kodaira dimension of moduli 
spaces of abelian surfaces'', 
{\it Compositio Math.} {\bf 72} (1989), no. 2, 121--163. 
\item{[Pu]} Puts, P.J., ``On some Fano-threefolds that are sections 
of Grassmannians'', {\it Nederl. Akad. Wetensch. Indag. Math.} 
{\bf 44}, (1982), no. 1, 77--90. 
\item{[Re]} Revoy, Ph.:
``Formes altern\'ees et puissances divis\'ees'', {\it S\'eminaire Dubreil}
1972/1973.
\item{[Schr]} Schreyer, F.-O., ``Geometry and algebra of prime Fano
3-folds of index 1 and genus 12'', preprint 1997.
\item{[Se]} Sekiguchi, T., ``On the cubics defining abelian varieties'', 
{\it J. Math. Soc. Japan},  {\bf 30}, (1978), no. 4, 703--721.
\item{[SR]} Semple, G., Roth, L., {\it Algebraic Geometry}, Chelsea 1937.
\item{[Si]} Silberger, A. J., 
``An elementary construction of the representations of 
${\rm SL}(2,\,{\rm GF}(q))$'', {\it Osaka J. Math.} {\bf 6}, (1969), 329--338. 
\item{[Tak]} Takeuchi, K: ``Some birational maps of Fano $3$-folds,'' 
{\it Compositio Math.}, {\bf 71}, (1989), no. 3, 265--283. 
\item{[Tan]} Tanaka, S., ``Construction and Classification of Irreducible
Representations of Special Linear Group of the Second Order over a Finite
Field,'' {\it Osaka J. Math.,} {\bf 4} (1967), 65--84.
\item{[Tre]} Tregub, S. L.: 
``Construction of a birational isomorphism of a three-dimensional cubic 
and a Fano variety of the first kind with $g=8$, 
connected with a normal rational curve of degree $4$,'' 
{\it Vestnik Moskov. Univ. Ser. I Mat. Mekh.}, (1985), no. 6, 99--101, {\bf 113}. 
(English translation: {\it Moscow Univ. Math. Bull.}, 
{\bf 40}, (1985), no. 6, 78--80.)
\item{[vdG]} van der Geer, G., ``On the geometry of a Siegel modular threefold,''
{\it Math. Ann.}, {\bf 260}, (1982), no. 3, 317--350. 
\item{[Ve]} V\'elu, J., ``Courbes elliptiques munies d'un sous-groupe
$\boldz/n\boldz\times\mu_n$'', {\it Bull. Soc. Math. France},
M\'emoire {\bf 57}, (1978) 5--152.
\item{[We]} Weil, A., ``Sur certains groupes d'op\'erateurs unitaires'', 
{\it Acta Math.}, {\bf 111}, (1964), 143--211. 
\end